\documentclass[a4paper,twoside]{amsart}
\usepackage{amsmath}
\usepackage{amsthm}
\usepackage{amssymb,amscd}
\input xy
\xyoption{all}
\numberwithin{equation}{section}
\theoremstyle{definition}
\newtheorem{dfn}{Definition}[section]
\newtheorem{example}[dfn]{Example}
\newtheorem{remark}[dfn]{Remark}
\theoremstyle{plain}
\newtheorem{thm}[dfn]{Theorem}
\newtheorem{prop}[dfn]{Proposition}
\newtheorem{cor}[dfn]{Corollary}
\newtheorem{lem}[dfn]{Lemma}
\newtheorem{claim}[dfn]{Claim}
\title[Locally finite simple weight modules over TGWAs]
{Locally finite simple weight modules over twisted
generalized Weyl algebras}
\author{Jonas T. Hartwig}
\address{Department of Mathematical Sciences, Division of Mathematics,
Chalmers University of Technology and
University of Gothenburg, Eklandagatan 86,
S-412 96 Gothenburg, Sweden}
\email{jonas.hartwig@math.chalmers.se}
\def\al{\alpha}
\def\be{\beta}
\def\ga{\gamma}
\def\ep{\varepsilon}

\def\la{\lambda}
\def\si{\sigma}
\def\th{\theta}
\def\ph{\varphi}
\def\fn{\mathfrak{n}}
\def\fm{\mathfrak{m}}
\def\fp{\mathfrak{p}}
\def\Z{\mathbb{Z}}
\def\Q{\mathbb{Q}}

\def\C{\mathbb{C}}
\def\N{\mathbb{N}}
\def\nn{\underline{n}}
\def\kk{\underline{k}}
\DeclareMathOperator{\supp}{supp}
\DeclareMathOperator{\rank}{rank}
\DeclareMathOperator{\spn}{Span}
\DeclareMathOperator{\sgn}{sgn}
\DeclareMathOperator{\Max}{Max}

\DeclareMathOperator{\GCD}{GCD}
\begin{document}
%\maketitle
%\begin{center}
%\begin{minipage}{75ex}
%\begin{footnotesize} \textsc{Abstract.} \,
%
%\end{footnotesize}
%\end{minipage}
%\end{center}
\begin{abstract}
We present methods and explicit formulas for describing
simple weight modules
over twisted generalized Weyl algebras.
When a certain commutative subalgebra is finitely generated over an algebraically
closed field we obtain a classification of a class of locally finite simple weight modules
as those induced from simple modules over a subalgebra
isomorphic to a tensor product of noncommutative tori.
As an application we describe simple weight modules over the quantized Weyl algebra.
\end{abstract}
\maketitle

%
%
%   I n t r o d u c t i o n
%
%
\section{Introduction}
Bavula defined in \cite{Ba}, \cite{B} the notion of a
generalized Weyl algebra (GWA) which is a class of algebras which
include $U(\mathfrak{sl}(2))$, $U_q(\mathfrak{sl}(2))$, down-up algebras,
and the Weyl algebra, as examples.
In addition to various ring theoretic properties,
the simple modules were also described for some
GWAs in \cite{B}.
In \cite{DGO} all simple and indecomposable weight
modules of GWAs of rank (or degree) one were
classified.

So called higher rank GWAs were defined in \cite{B} and
%as tensor products of
%rank one GWAs. This has some consequences on the
%side of representations.
in \cite{B2} the authors studied indecomposable weight modules
over certain higher rank GWAs.

In \cite{MT}, with the goal to enrich the
representation theory in the higher rank case,
the authors defined the twisted generalized
Weyl algebras (TGWA). This is a class of algebras which
include all higher rank GWAs (if a certain
subring $R$ has no zero divisors) and also many
algebras which can be viewed as twisted
tensor products of rank one GWAs, for example certain
Mickelsson step algebras and extended Orthogonal
Gelfand-Zetlin algebras \cite{MPT}.
Under a technical assumption on the algebra
formulated using a biserial graph,
some torsion-free simple weight modules were
described in \cite{MT}.
Simple graded weight modules were studied in
\cite{MPT} using an analogue of the Shapovalov
form.

In this paper we describe a more general class
of locally finite simple weight modules
over TGWAs using the well-known technique of considering
the maximal graded subalgebra which preserves
the weight spaces. It is known that under quite
general assumptions (see Theorem 18 in \cite{DFO})
any simple weight module over a TGWA is a unique quotient of
a module which is induced from a simple module over
this subalgebra.
Our main results are the description of
this subalgebra under various assumptions
(Theorem \ref{thm:varphibmstruct1}
and Theorem \ref{thm:Bisotensor})
and the explicit formulas (Theorem \ref{thm:action})
of the associated module of the TGWA.
In contrast to \cite{MT}, we do not assume that
the orbits are torsion-free and we allow the modules
to have some inner breaks, as long as they do not have
any so called \emph{proper} inner breaks (see Definition \ref{dfn:NPIB}).
The weight spaces will not in general be one-dimensional in
our case, which was the case in \cite{MT}, \cite{MPT}.

Moreover, as an application we classify the simple weight modules
without proper inner breaks over a quantized Weyl algebra
of rank two (Theorem \ref{thm:QWAcl}).

The paper is organized as follows.
In Section \ref{sec:definitions} the definitions of twisted
generalized Weyl constructions and algebras are
given together with some examples.
Weight modules and the subalgebra $B(\omega)$
are defined.

In Section \ref{sec:preliminaries} we first prove
some simple facts and then define the class of
simple weight modules with no proper inner breaks.
We also show that this class properly contains all the
modules studied in \cite{MT}.

Section \ref{sec:subalg} is devoted to the description
of the subalgebra $B(\omega)$.
When the ground field is algebraically closed and a certain subalgebra
$R$ is finitely generated,
we show that it is
isomorphic to a tensor product of noncommutative tori
for which the finite-dimensional irreducible representations
are easy to describe.

In Section \ref{sec:explicit} we specify a basis and give
explicit formulas for the irreducible quotient of the
induced module.
%lift the action
%from the subalgebra to the TGWA by specifying a
%basis and giving explicit formulas
%for the action.

Finally, in Section \ref{sec:qweyl} we consider
as an example the quantized Weyl algebra and
determine certain important subsets of $\Z^n$
related to $B(\omega)$ and the support of modules
as solutions to some systems of equations.
In the rank two case we describe all simple weight modules
with finite-dimensional weight spaces and no
proper inner breaks.

%
%
%   D e f i n i t i o n s
%
%
\section{Definitions}\label{sec:definitions}

%
% The TGWC and TGWA
%
\subsection{The TGWC and TGWA}
Fix a positive integer $n$ and set $\nn=\{1,2,\ldots,n\}$.
Let $K$ be a field,
and let $R$ be a commutative unital $K$-algebra,
$\boldsymbol{\si}=(\si_1,\ldots,\si_n)$ be an $n$-tuple of pairwise
commuting $K$-automorphisms of $R$,
$\boldsymbol{\mu}=(\mu_{ij})_{i,j\in\nn}$ be a matrix with
entries from $K^*:=K\backslash\{0\}$
and $\boldsymbol{t}=(t_1,\ldots,t_n)$ be an $n$-tuple of
nonzero elements from $R$.
The \emph{twisted generalized Weyl construction} (TGWC)
$A'$ obtained from
the data $(R,\boldsymbol{\si},\boldsymbol{t},\boldsymbol{\mu})$
%(\s_i)_{i\in\nn},(t_i)_{i\in\nn},(\m_{i,j})_{i,j\in\nn})$
is the unital $K$-algebra generated over $R$ by $X_i, Y_i, (i\in\nn)$
with the relations
\begin{align}
X_ir&=\si_i(r)X_i,&Y_ir&=\si_i^{-1}(r)Y_i, &\text{for }r\in R, i\in\nn,
\label{eq:rel1}\\
Y_iX_i&=t_i, &X_iY_i&=\si_i(t_i), &\text{for }i\in\nn,\\
X_iY_j&=\mu_{ij}Y_jX_i, &&&\text{for }i,j\in\nn, i\neq j.
\label{eq:rel3}
\end{align}
From the relations (\ref{eq:rel1})--(\ref{eq:rel3}) follows that $A'$ carries a
$\Z^n$-gradation $\{A'_g\}_{g\in\Z^n}$ which is uniquely defined by requiring
\[\deg X_i=e_i, \quad\deg Y_i=-e_i,\quad \deg r=0,\quad\text{for }i\in\nn,
r\in R,\]
where $e_i=(0,\ldots,\overset{i}{1},\ldots,0)$.
The \emph{twisted generalized Weyl algebra} (TGWA)
$A=A(R,\boldsymbol{\si},\boldsymbol{t},\boldsymbol{\mu})$
of \emph{rank} $n$ is defined to be $A'/I$, where $I$
is the sum of all graded two-sided ideals of $A'$ intersecting $R$
trivially. Since $I$ is graded, $A$ inherits a $\Z^n$-gradation
$\{A_g\}_{g\in\Z^n}$ from $A'$.

Note that from relations (\ref{eq:rel1})--(\ref{eq:rel3}) follows the
identity
\begin{equation}\label{eq:xixj}
X_iX_jt_i=X_jX_i\mu_{ji}\si_j^{-1}(t_i)
\end{equation}
which holds for $i,j\in\nn, i\neq j$.
Multiplying (\ref{eq:xixj}) from the left by $\mu_{ij}Y_j$ we obtain
\begin{equation}\label{eq:xititj}
X_i\big(t_it_j-\mu_{ij}\mu_{ji}\si_i^{-1}(t_j)\si_j^{-1}(t_i)\big)=0
\end{equation}
for $i,j\in\nn, i\neq j$.
One can show that the algebra $A'$, hence $A$, is nontrivial if one assumes that
$t_it_j=\mu_{ij}\mu_{ji}\si_i^{-1}(t_j)\si_j^{-1}(t_i)$
for $i,j\in\nn, i\neq j$. Analogous identities exist for $Y_i$.

\subsection{Examples}\label{sec:examples}
Some of the first motivating examples of a \emph{generalized Weyl algebra}
(GWA), i.e. a TGWC of rank $1$, are $U(\mathfrak{sl}(2))$,
$U_q(\mathfrak{sl}(2))$ and of course the Weyl algbra $A_1$. We
refer to \cite{B} for details.

We give some examples of TGWAs of higher rank.
\subsubsection{Quantized Weyl algebras}
Let $\Lambda=(\la_{ij})$ be an $n\times n$ matrix with
nonzero complex entries such that $\la_{ij}=\la_{ji}^{-1}$.
Let $\bar q=(q_1,\ldots,q_n)$ be an $n$-tuple of elements
of $\C\backslash\{0,1\}$. The $n$:th quantized Weyl algebra
$A_n^{\bar q,\Lambda}$ is the $\C$-algebra with generators
$x_i,y_i$, $1\le i\le n$, and relations
\begin{align}
x_ix_j&=q_i\la_{ij}x_jx_i, & y_iy_j&=\la_{ij}y_jy_i,
\label{eq:qweylrel1}\\
x_iy_j&=\la_{ji}y_jx_i, & x_jy_i&=q_i\la_{ij}y_ix_j,
\label{eq:qweylrel2}
\end{align}
\begin{equation}
\label{eq:qweylrel3}
x_iy_i-q_iy_ix_i=1+\sum_{k=1}^{i-1}(q_k-1)y_ix_i,
\end{equation}
for $1\le i<j\le n$.
Let $R=\C[t_1,\ldots, t_n]$ be the polynomial algebra in
$n$ variables and $\si_i$ the $\C$-algebra
automorphisms defined by
\begin{equation}\label{eq:qweylsigmadef}
\si_i(t_j)=
\begin{cases}
t_j, & j<i, \\
1+q_it_i+\sum_{k=1}^{i-1}(q_k-1)t_k, & j=i, \\
q_it_j, & j>i.
\end{cases}
\end{equation}
One can check that the $\si_i$ commute.
Let $\boldsymbol{\mu}=(\mu_{ij})_{i,j\in\nn}$
be defined by $\mu_{ij}=\la_{ji}$ and $\mu_{ji}=q_i\la_{ij}$
for $i<j$. Let also
$\boldsymbol{\si}=(\si_1,\ldots,\si_n)$ and
$\boldsymbol{t}=(t_1,\ldots,t_n)$.
One can show that the maximal graded ideal of
the TGWC
$A'(R,\boldsymbol{\si},\boldsymbol{t},\boldsymbol{\mu})$
is generated by the elements
\[X_iX_j-q_i\la_{ij}X_jX_i,\; Y_iY_j-\la_{ij}Y_jY_i,\quad
1\le i<j\le n.\]
Thus $A_n^{\bar q,\Lambda}$ is isomorphic
to the TGWA $A(R,\boldsymbol{\si},\boldsymbol{t},\boldsymbol{\mu})$
via $x_i\mapsto X_i$, $y_i\mapsto Y_i$.

\subsubsection{$Q_{ij}$-CCR}
Let $(Q_{ij})_{i,j=1}^d$ be an $d\times d$ matrix with
complex entries such that $Q_{ij}=Q_{ji}^{-1}$
if $i\neq j$ and
$A_d$ be the algebra generated by elements $a_i,a_i^*$,
$1\le i\le d$ and relations
\begin{align*}
a_i^*a_i-Q_{ii}a_ia_i^*&=1, &a_i^*a_j&=Q_{ij}a_ja_i^*, \\
a_ia_j&=Q_{ji}a_ja_i, & a_i^*a_j^*&=Q_{ij}a_j^*a_i^*,
\end{align*}
where $1\le i,j\le d$ and $i\neq j$.
Let $R=\C[t_1,\ldots, t_d]$ and define
the autormorphisms $\si_i$ of $R$ by $\si_i(t_j)=t_j$
if $i\neq j$ and $\si_i(t_i)=1+Q_{ii}t_i$.
Let $\mu_{ij}=Q_{ji}$ for all $i,j$. Then
$A_d$ is isomorphic to the TGWA
$A(R,(\si_1,\ldots,\si_n),(t_1,\ldots,t_n),
\boldsymbol{\mu})$.

\subsubsection{Mickelsson and OGZ algebras}
In both the above examples the generators $X_i$ and $X_j$
commute up to a multiple of the ground field. This need
not be the case as shown in \cite{MPT}, where it was
shown that Mickelsson step algebras and extended
orthogonal Gelfand-Zetlin algebras are TGWAs.

%
%   W e i g h t   m o d u l e s
%
\subsection{Weight modules}
Let $A$ be a TGWC or a TGWA.
Let $\Max(R)$ denote the set of all maximal ideals in $R$. A module $M$ over $A$
is called a \emph{weight module} if
\[M=\oplus_{\fm\in\Max(R)} M_\fm,\]
where
\[M_\fm=\{v\in M\;|\; \fm v = 0\}.\]
The \emph{support}, $\supp(M)$, of $M$ is the set of all $\fm\in\Max(R)$ such
that $M_\fm\neq 0$.
A weight module is \emph{locally finite} if all the weight spaces
$M_\fm$, $\fm\in\supp(M)$, are finite-dimensional over the ground field $K$.

Since the $\si_i$ are pairwise commuting, the free abelian group $\Z^n$
acts on $R$ as a group of $K$-algebra automorphisms by
\begin{equation}\label{eq:Znaction}
g(r)=\si_1^{g_1}\si_2^{g_2}\ldots\si_n^{g_n}(r)
\end{equation}
for $g=(g_1,\ldots,g_n)\in \Z^n$ and $r\in R$.
Then $\Z^n$ also acts naturally on $\Max(R)$ by
$g(\fm)=\{g(r)\;|\;r\in \fm\}$.
Note that
\begin{equation} \label{eq:XiMm}
X_iM_\fm\subseteq M_{\si_i(\fm)}\quad\text{and}\quad
Y_iM_\fm\subseteq M_{\si_i^{-1}(\fm)}
\end{equation}
for any $\fm\in\Max(R)$.
If $a\in A$ is homogenous of degree $g\in\Z^n$, then by using
(\ref{eq:rel1}) and (\ref{eq:XiMm}) repeatedly one obtains
the very useful identities
\begin{equation}\label{eq:araMm1}
a\cdot r=g(r)\cdot a,\quad r\cdot a=a\cdot (-g)(r),
\end{equation}
for $r\in R$ and
\begin{equation}\label{eq:araMm2}
aM_\fm\subseteq M_{g(\fm)}
\end{equation}
for $\fm\in\Max(R)$.

%
%
% S u b a l g e b r a s   p r e s e r v i n g   w e i g h t s p a c e s
%
%
\subsection{Subalgebras leaving the weight spaces invariant}
\label{sec:subalgleaving}
Let $\omega\subseteq\Max(R)$ be an orbit under the action of $\Z^n$
on $\Max(R)$ defined in (\ref{eq:Znaction}).
Let
\begin{equation}\label{eq:Znga}
\Z^n_\omega=\Z^n_\fm=\{g\in \Z^n\; |\; g(\fm)=\fm\}
\end{equation}
where $\fm$ is some point in $\omega$.
Since $\Z^n$ is abelian, $\Z^n_\omega$ does not depend on the choice of $\fm$ from
$\omega$.
Define
\begin{equation}\label{eq:Bga}
B(\omega)=\oplus_{g\in \Z^n_\omega}A_g.
\end{equation}
Since $A$ is $\Z^n$-graded and since $\Z^n_\omega$ is a subgroup of $\Z^n$,
$B(\omega)$ is a subalgebra of $A$ and $R=A_0\subseteq B(\omega)$.
Let $\fm\in\omega$ and suppose that $M$ is a simple weight $A$-module
with $\fm\in\supp(M)$.
Since $M$ is simple we have $\supp(M)\subseteq \omega$.
Using (\ref{eq:araMm2}) it follows that $B(\omega)M_\fm\subseteq M_\fm$
and by definition $M_\fm$ is annihilated by $\fm$ hence also by
the two-sided ideal $(\fm)$ in $B(\omega)$ generated by $\fm$.
Thus $M_\fm$ is naturally a module over the algebra
\begin{equation}\label{eq:Bfm}
B_\fm:=B(\omega)/(\fm).
\end{equation}
By Proposition 7.2 in \cite{MPT} (see also Theorem 18 in \cite{DFO} for
a general result),
$M_\fm$ is a simple $B_\fm$-module, and any simple $B_\fm$-module
occurs as a weight space in a simple weight $A$-module.
Moreover, two simple weight $A$-modules $M,N$ are isomorphic if
and only if $M_\fm$ and $N_\fm$ are isomorphic as $B_\fm$-modules.
Therefore we are led
to study the algebra $B_\fm$ and simple modules over it.

%
%
%   P r e l i m i n a r i e s
%
%
\section{Preliminaries}\label{sec:preliminaries}
%
%   R e d u c e d   w o r d s
%
\subsection{Reduced words}
Let $L=\{X_i\}_{i\in\nn}\cup\{Y_i\}_{i\in\nn}$.
By a \emph{word} $(a; Z_1,\ldots, Z_k)$ in $A$ we will mean an element
$a$ in $A$ which is a product of elements from the set $L$,
together with a fixed tuple $(Z_1,\ldots,Z_k)$ of elements
from $L$ such that $a=Z_1\cdot \ldots\cdot Z_k$. When referring to a
word we will often
write $a=Z_1\ldots Z_k\in A$ to denote the word $(a;Z_1,\ldots,Z_k)$
or just write $a\in A$, suppressing
the fixed representation of $a$ as a product of elements from $L$.

Set $X_i^*=Y_i$ and $Y_i^*=X_i$. For a word
$a=Z_1\ldots Z_k\in A$ we define
\[a^*:=Z_k^*\cdot Z_{k-1}^*\cdot\ldots\cdot Z_1^*.\]

In the special case when $\mu_{ij}=\mu_{ji}$ for all $i,j$ then
by (\ref{eq:rel1})-(\ref{eq:rel3}) there is an
anti-involution $*$
on $A'$ defined by $X_i^*=Y_i$, and $r^*=r$ for $r\in R$. Since $I^*=I$
this anti-involution carries over to $A$.

\begin{dfn}
A word $Z_1\ldots Z_k$ will be called \emph{reduced}
if
\[Z_i\neq Z_j^* \text{ for } i,j\in\kk\]
and
\[Z_i\in\{X_r\}_{r\in\nn} \Longrightarrow Z_j\in\{X_r\}_{r\in\nn}
\;\forall j\ge i.\]
\end{dfn}

For example $Y_1Y_2Y_1X_3$ is reduced whereas $Y_1Y_2X_1$ and $Y_1X_2Y_3$
are not. The following Lemma and Corollary explains the importance of the
reduced words.

\begin{lem} \label{lem:puresep}
Any word $b$ in $A$ can be written $b=a\cdot r=r'\cdot a$, where
$a$ is a reduced word, and $r,r'\in R$.
\end{lem}
\begin{proof} If $a$ and $r$ has been found we can take
$r'=(\deg a)(r)$, according to (\ref{eq:araMm1}). Thus we
concentrate on finding $a$ and $r$. Let $b=Z_1\ldots Z_k$ be an
arbitrary word in $A$. We prove the statement by induction on $k$.
If $k=1$, then $b$ is necessarily reduced so take $a=b, r=1$. When
$k>1$, use the induction hypothesis to write
\[Z_1\ldots Z_{k-1}=Y_{i_1}\ldots Y_{i_l}X_{j_1}\ldots X_{j_m}\cdot r',\]
where $1\le i_u,j_v\le n$ and $i_u\neq j_v$ for any $u,v$.
Consider first the case when $Z_k=Y_j$ for some $j\in\nn$.
Then
\[Z_1\ldots Z_k=Y_{i_1}\ldots Y_{i_l}X_{j_1}\ldots X_{j_m}Y_j\cdot \si_j(r').\]
If $j_v\neq j$ for $v=1,\ldots,m$ we are done because
using relation (\ref{eq:rel3}) repeatedly we obtain,
\[Z_1\ldots Z_k=Y_{i_1}\ldots Y_{i_l} Y_j X_{j_1}\ldots X_{j_m}\cdot
\mu\si_j(r')\]
for some $\mu\in K^*$.
Otherwise, let $v\in\{1,\ldots,m\}$ be maximal such that $j_v=j$. Then
\begin{align*}
Z_1\ldots Z_k&=Y_{i_1}\ldots Y_{i_l} X_{j_1}\ldots
X_{j_v} Y_j X_{j_{v+1}}\ldots X_{j_m} \mu \si_j(r')=\\
&=Y_{i_1}\ldots Y_{i_l} X_{j_1}\ldots X_{j_{v-1}}X_{j_{v+1}}\ldots X_{j_m}
 w(t_j)\mu \si_j(r')
\end{align*}
for some $\mu\in K^*$ and some $w\in W$. It remains to consider
the case $Z_k=X_j$ for some $j\in\nn$. But using that
\[Y_{i_1}\ldots Y_{i_l}X_{j_1}\ldots X_{j_m}=X_{j_1}\ldots X_{j_m}Y_{i_1}\ldots
Y_{i_l}\mu\]
for some $\mu\in K^*$, it is clear that this case is analogous.
\end{proof}

\begin{cor} \label{cor:genAd}
Each $A_g$, $g\in W$, is generated as a right (and also as a
left) $R$-module by the reduced words of degree $g$.
\end{cor}

\begin{lem} \label{lem:primestar}
Suppose $*$ defines an anti-involution on $A$.
Let $\fp$ be a prime ideal of $R$.
Let $g\in\Z^n$ and let $a\in A_g$.
If $ba\notin \fp$ for some $b\in A_{-g}$ then $a^*a\notin\fp$.
\end{lem}
\begin{proof} Since $\fp$ is prime, and $ba\in R$ we have
\[\fp\not\ni (ba)^2=(ba)^*ba=a^*b^*ba=a^*a\cdot(-\deg a)(b^*b)\]
so in particular $a^*a\notin\fp$.
\end{proof}
\begin{remark}\label{rem:primestar}
If we assume $a$ and $b$ to be words in the formulation of 
Lemma \ref{lem:primestar}, one can easily show that the
statement remains true without the restriction on $*$ to
be an anti-involution.
\end{remark}

%
%   I n n e r   b r e a k s
%
\subsection{Inner breaks and canonical modules}
Let $A$ be a TGWC or a TGWA and let $M$ be a simple weight module over $A$.
In \cite{MT} Remark 1 it was noted that the problem of describing simple
weight modules over a TGWC is wild in general. This is a motivation for
restricting attention to some subclass which has nice properties.
In \cite{MT} the following definition was made.
\begin{dfn}
The support of $M$ has \emph{no inner breaks} if for all $\fm\in\supp(M)$,
\begin{align*}
t_i\in\fm\Longrightarrow\si_i(\fm)\notin\supp(M),\text{ and}\\
\si_i(t_i)\in\fm\Longrightarrow\si_i^{-1}(\fm)\notin\supp(M).
\end{align*}
\end{dfn}
We introduce the following property.
\begin{dfn}\label{dfn:NPIB}
We say that $M$ has \emph{no proper inner breaks} if
for any $\fm\in\supp(M)$ and any word $a$ with $aM_\fm\neq 0$
we have $a^*a\notin\fm$.
\end{dfn}

Observe that whether or not $a^*a\in\fm$ for a word $a$ does not
depend on the particular representation of $a$ as a product
of generators.
Note also that to prove that
a simple weight module $M$ has no proper inner breaks,
it is sufficient
to find for any $\fm\in\supp(M)$ and any word $a$ with
$aM_\fm\neq 0$ a word $b\in A$
of degree $-\deg a$ such that $ba\notin\fm$
because then $a^*a\notin\fm$ automatically
by Remark \ref{rem:primestar}.
In fact one can show that a simple weight module $M$
has no proper inner breaks if (and only if) there exists
an $\fm\in\supp(M)$ such that for any reduced word
$a\in A$ with $aM_\fm\neq 0$ and $aM_\fm\subseteq M_\fm$
there is a word $b$ of degree $-\deg a$ such
that $ba\notin\fm$. However we will not use this result.

The choice of terminology in Definition \ref{dfn:NPIB} is motivated
by the following proposition.
\begin{prop}
If $M$ has no inner breaks, then $M$ has no proper
inner breaks either.
\end{prop}
\begin{proof}
Let $\fm\in\supp(M)$ and
$a=Z_1\ldots Z_k\in A$ be a word such that $aM_\fm\neq 0$.
Thus $Z_i\ldots Z_k M_\fm\neq 0$ for $i=1,\ldots,k+1$
so (\ref{eq:araMm2}) implies that
\[(\deg Z_i\ldots Z_k)(\fm)\in\supp(M).\]
If $M$ has no inner breaks, it follows that
$Z_i^*Z_i\notin (\deg Z_{i+1}\ldots Z_k)(\fm)$
for $i=1,\ldots, k$.
Now using (\ref{eq:araMm1}),
\begin{align}
a^*a=Z_k^*\ldots Z_1^*Z_1\ldots Z_k&=
Z_k^*\ldots Z_2^*Z_2\ldots Z_k
(-\deg Z_2\ldots Z_k)(Z_1^*Z_1)=\nonumber\\
&=\ldots=
\prod_{i=1}^k (-\deg Z_{i+1}\ldots Z_k)(Z_i^*Z_i)\notin\fm.
\label{eq:astara}
\end{align}
Thus $M$ has no proper inner breaks.
\end{proof}

In \cite{MT}, a simple weight module $M$ was defined to be
\emph{canonical} if for any $\fm,\fn\in\supp(M)$ there is
an automorphism $\si$ of $R$ of the form
\[\si=\si_{i_1}^{\ep_1}\cdot\ldots\cdot\si_{i_k}^{\ep_k},
\quad\ep_j=\pm 1\text{ and }1\le i_j\le n,\text{ for }j=1,\ldots,k,\]
such that $\si(\fm)=\fn$ and
such that for each $j=1,\ldots,k$,
\begin{equation}\label{eq:can1}
t_{i_j}\notin \si_{i_{j+1}}^{\ep_{j+1}}\ldots\si_{i_k}^{\ep_k}(\fm)
\quad\text{if }\ep_j=1,\text{ and}
\end{equation}
\begin{equation}\label{eq:can2}
\si_{i_j}(t_{i_j})\notin
\si_{i_{j+1}}^{\ep_{j+1}}\ldots\si_{i_k}^{\ep_k}(\fm)
\quad\text{if }\ep_j=-1.
\end{equation}
This definition can be reformulated as follows.
\begin{prop} $M$ is canonical iff for any $\fm,\fn\in\supp(M)$
there is a word $a\in A$ such that
$aM_\fm\subseteq M_\fn$ and $a^*a\notin\fm$.
\end{prop}
\begin{proof} Suppose $M$ is canonical, and let
$\fm,\fn\in\supp(M)$. Let $\si$ be as in the definition
of canonical module. Define $a=Z_1\ldots Z_k$
where $Z_j=X_{i_j}$ if $\ep_j=1$ and $Z_j=Y_{i_j}$ otherwise.
Using (\ref{eq:araMm2}) we see that $aM_\fm\subseteq M_\fn$.
Also, (\ref{eq:can1}) and (\ref{eq:can2}) translates into
\[Z_j^*Z_j\notin (\deg Z_{j+1}\ldots Z_k)(\fm)\]
for $j=1,\ldots,k$. Using the calculation (\ref{eq:astara})
and that $\fm$ is prime we deduce that $a^*a\notin\fm$.

Conversely, given a word $a=Z_1\ldots Z_k\in A$
with $aM_\fm\subseteq M_\fn$ and $a^*a\notin\fm$,
we define $\ep_i=1$ if $Z_i=X_i$ and $\ep_i=-1$ otherwise.
Then from $a^*a\notin\fm$ follows that
$\si:=\si_{i_1}^{\ep_1}\cdot\ldots\cdot\si_{i_k}^{\ep_k}$
satisfies (\ref{eq:can1}) and (\ref{eq:can2})
by the same reasoning as above.
\end{proof}
\begin{cor}\label{cor:noPIBcan}
If $M$ has no proper inner breaks, then $M$ is canonical.
\end{cor}
\begin{proof} We only need to note that
since $M$ is a simple weight module there
is for each $\fm,\fn\in\supp(M)$ a word $a$ such that
$0\neq aM_\fm\subseteq M_\fn$.
\end{proof}

Under the assumptions in \cite{MT} any canonical module
has no inner breaks (see \cite{MT}, Proposition 1).
However we have the following example of a TGWA $A$
and a simple weight module $M$ over $A$ which has
no proper inner breaks, and thus is canonical by
Corollary \ref{cor:noPIBcan}, but nonetheless
has an inner break.
\begin{example} Let
$R=\C[t_1,t_2]$ and define the $\C$-algebra automorphisms
$\si_1$ and $\si_2$ of $R$ by
$\si_i(t_j)=-t_j$ for $i,j=1,2$. Let
$\boldsymbol{\mu}=
[\begin{smallmatrix}0&1\\1&0\end{smallmatrix}]$.
Let $A'=A'(R, \boldsymbol{t}, \boldsymbol{\si}, \boldsymbol{\mu})$
be the associated TGWC, where $\boldsymbol{t}=(t_1,t_2),
\boldsymbol{\si}=(\si_1,\si_2)$.
Then one can check that
$I=\langle X_1X_2+X_2X_1, Y_1Y_2+Y_2Y_1\rangle$.
Let $M$ be a vector space over $\C$ with basis
$\{v,w\}$ and define an $A'$-module structure
on $M$ by letting $X_1M=Y_1M=0$ and
\begin{align*}
X_2v&=w,&X_2w&=v,\\
Y_2v&=w,&Y_2w&=-v.
\end{align*}
It is easy to check that the required relations
are satisfied and that $IM=0$, hence $M$ becomes
an $A$-module. Let
$\fm=(t_1,t_2+1)$ and $\fn=(t_1,t_2-1)$. Then
\[M=M_\fm\oplus M_\fn, \quad\text{where } M_\fm=\C v, M_\fn=\C w\]
so $M$ is a weight module. Any proper nonzero submodule of $M$
would also be a weight module by standard results. That no
such submodule can exist is easy to check, so
$M$ is simple. One checks that $M$ has no proper
inner breaks. But $t_1\in\fm$ and
$\si_1(\fm)=\fn\in\supp(M)$ so $\fm$ is an
inner break.
\end{example}

%
%
% The weight space preserving subalgebra and its representations
%
%
\section{The weight space preserving subalgebra and its irreducible representations}
\label{sec:subalg}
In this section, let $A$ be a TGWC,
$\fm\in\Max(R)$ and let $\omega$ be the $\Z^n$-orbit of
$\fm$. Recall the set $\Z^n_\omega$ defined in (\ref{eq:Znga}).
Define the following subsets of $\Z^n$:
\begin{equation}\label{eq:tildeGmandGm}
\tilde G_\fm=\{g\in \Z^n\; |\; a^*a\notin\fm\text{ for some word }a\in A_g\}
\quad \text{and}\quad G_\fm=\tilde G_\fm\cap \Z^n_\omega.
\end{equation}
Let also $\ph_\fm:A\to A/(\fm)$ denote the canonical projection,
where $(\fm)$ is the two-sided ideal in $A$ generated by
$\fm$, and let
$R_\fm=R/\fm$ be the residue field of $R$ at $\fm$.

\begin{lem}\label{lem:dim}
Let $g\in\tilde G_\fm$. Then
\begin{equation}\label{eq:lemdim}
\ph_\fm(A_g)=R_\fm\cdot \ph_\fm(a)=\ph_\fm(a)\cdot R_\fm
\end{equation}
for any word $a\in A_g$ with $a^*a\notin\fm$.
\end{lem}
\begin{proof}
Let $b\in A_g$ be any element and
$a\in A_g$ a word such that $a^*a\notin\fm$,
We must show that there is an $r\in R$ such that
$\ph_\fm(b)=\ph_\fm(r)\ph_\fm(a)$.
Since $a^*a\notin\fm$ and $\fm$ is maximal,
$1-r_1a^*a\in\fm$ for some $r_1\in R$.
Set $r=br_1a^*$. Then $r\in R$ and
\[b-ra=b(1-r_1a^*a)\in (\fm).\]
The last equality in (\ref{eq:lemdim}) is immediate using (\ref{eq:araMm1}).
\end{proof}

The following result was proved in \cite{MT} Lemma 8 for
simple weight modules with so called regular support which
in particular means that they have no inner breaks.
It is still true in the more general situation when
$M$ has no proper inner breaks. Recall the ideal $I$ from the
definition of a $TGWA$.
\begin{prop} \label{prop:IM=0}
Suppose $A$ is a TGWC. If $M$ is a simple weight $A$-module
with no proper inner breaks, then $IM=0$. Hence $M$ is naturally a module over
the associated TGWA $A/I$.
\end{prop}
\begin{proof} Since $I$ is graded and $M$ is a weight modules,
it is enough to show that $(I\cap A_g)M_\fm=0$ for any $g\in\Z^n$
and any $\fm\in\supp(M)$. Assume that $a\in I\cap A_g$ and $av\neq 0$
for some $v\in M_\fm$. Then $a_1v\neq 0$ for some word $a_1$ in $a$.
Since $M$ has no proper inner breaks, $a_1^*a_1\notin\fm$ so by
Lemma \ref{lem:dim} there is an $r\in R$ such that $av=a_1rv$.
Thus $0\neq a_1^*a_1rv=a_1^*av$ which implies that $a_1^*a\in R\backslash\{0\}$.
This contradicts that $a\in I$.
\end{proof}

We fix now for each $g\in\tilde G_\fm$ a word $a_g\in A_g$ such that
$a_g^*a_g\notin\fm$. For $g=0$ we choose $a_g=1$.

\begin{lem}\label{lem:NBabc}
For any $g\in \tilde G_\fm, h\in G_\fm$ we have\\
a) $(a_ga_h^*)^*a_ga_h^*\notin\fm$ so in particular $g-h\in \tilde G_\fm$ and $G_\fm$ is
a subgroup of $\Z^n_\omega$,\\
b) $\ph_\fm(A_g)\ph_\fm(A_h)=\ph_\fm(A_gA_h)=\ph_\fm(A_{g+h})$,\\
c) $A_{g+h}M_\fm=A_gM_\fm$.
\end{lem}
\begin{proof}
a) We have
\begin{equation}\label{eq:abc1}
(a_ga_h^*)^*a_ga_h^*=a_ha_g^*a_ga_h^*=a_ha_h^*h(a_g^*a_g).
\end{equation}
Now $a_g^*a_g\notin\fm$ so $h(a_g^*a_g)\notin h(\fm)=\fm$. And
\[\fm\not\ni (a_h^*a_h)^2=a_h^*(a_ha_h^*)a_h=a_h^*a_h\cdot (-h)(a_ha_h^*)\]
so $a_ha_h^*\notin h(\fm)=\fm$. Since $\fm$ is maximal the right
hand side of (\ref{eq:abc1}) does not belong to $\fm$.
Since $\deg(a_ga_h^*)=g-h$ we obtain
$g-h\in \tilde G_\fm$. If in addition $g\in G_\fm$ then $g-h\in \Z^n_\omega$ also
since $\Z^n_\omega$ is a group. Thus $g-h\in G_\fm$ so
$G_\fm$ is a subgroup of $\Z^n_\omega$.

b) Since $\ph_\fm$ is a homomorphism, the first equality holds. By part a), $-h\in G_\fm$ so by
part a) again, $(a_ga_{-h}^*)^*a_ga_{-h}^*\notin\fm$. Hence
by Lemma \ref{lem:dim}, we have
\[\ph_\fm(A_{g+h})=R_\fm\cdot\ph_\fm(a_ga_{-h}^*)\subseteq \ph_\fm(A_gA_h).\]
The reverse inclusion holds since $\{A_g\}_{g\in\Z^n}$ is a gradation of $A$.

c) By part a), $g+h=g-(-h)\in \tilde G_\fm$. Thus by part b),
\[A_{g+h}M_\fm=\ph_\fm(A_{g+h})M_\fm=
\ph_\fm(A_gA_h)M_\fm=
A_gA_hM_\fm\subseteq A_gM_{h(\fm)}=A_g M_\fm.\]
By part a), the same calculation holds if we replace $g$ by $g+h$ and
and $h$ by $-h$, which gives the opposite inclusion.
\end{proof}

\begin{lem}\label{lem:gmtilde}
Let $g\in\Z^n\backslash\tilde G_\fm$. Then $A_gM_\fm=0$
for any simple weight module $M$ over $A$
with no proper inner breaks.
\end{lem}
\begin{proof}
Let $a\in A_g$ be any word. Then
$a^*a\in\fm$ and hence if $M$ is a simple weight module over $A$ with
no proper inner breaks, $aM_\fm=0$. Since the words generate $A_g$ as
a left $R$-module, it follows that $A_gM_\fm=0$.
\end{proof}

%
%
%  S t r u c t u r e   o f   $B_\fm$  -  G e n e r a l   c a s e
%
%
\subsection{General case}\label{sec:generalcase}
Recall that $(\fm)$ denotes the two-sided ideal in $A$ generated by $\fm$.
Since $(\fm)$ is a graded ideal in $A$, there is an induced $\Z^n$-gradation
of the quotient $A/(\fm)$ and $\ph_\fm(A_g)=(A/(\fm))_g$.
Corresponding to the decomposition $\Z^n_\omega$ into the subset $G_\fm$ and
its complement are two $K$-subspaces of the algebra
$B_\fm=B(\omega)/\big(B(\omega)\cap (\fm)\big)$
which will be denoted by $B_\fm^{(1)}$ and $B_\fm^{(0)}$ respectively.
In other words, $B_\fm=B_\fm^{(1)}\oplus B_\fm^{(0)}$, where
\[ B_\fm^{(1)}=\bigoplus_{g\in G_\fm} (A/(\fm))_g
\quad\text{and}\quad
B_\fm^{(0)}=\bigoplus_{g\in\Z^n_\omega\backslash G_\fm} (A/(\fm))_g.\]

By Lemma \ref{lem:NBabc}a), $G_\fm$ is a subgroup of the free abelian
group $\Z^n$, hence is free abelian itself of rank $k\le n$.
Let $s_1,\ldots, s_k$ denote a basis for $G_\fm$ over $\Z$ and let
$b_i=\ph_\fm(a_{s_i})$ for $i=1,\ldots,k$.
Note also that $R_\fm$ is an extension field of $K$ and that $\Z^n_\omega$
acts naturally on $R_\fm$ as a group of $K$-automorphisms. Let
$\{\rho_j\}_{j\in J}$ be a basis for $R_\fm$ over $K$.

\begin{thm} \label{thm:varphibmstruct1}
a) $B_\fm^{(0)}M_\fm=0$ for any simple weight module $M$ over $A$
with no proper inner breaks, and\\
b) the $b_i$ are invertible and as a $K$-linear space,
$B_\fm^{(1)}$ has a basis
\begin{equation}\label{eq:bm1basis}
\{\rho_j b_1^{l_1}\ldots b_k^{l_k} \; |\; j\in J\text{ and } l_i\in\Z\text{ for } 1\le i\le k\}
\end{equation}
and the following commutation relations hold
\begin{equation} \label{eq:Brel1}
b_i\lambda=s_i(\lambda) b_i, \quad i=1,\ldots,k, \la\in R_\fm,
\end{equation}
\begin{equation} \label{eq:Brel2}
b_ib_j=\lambda_{ij}b_jb_i,\quad i,j=1,\ldots,k
\end{equation}
for some nonzero $\lambda_{ij}\in R_\fm$.
\end{thm}
\begin{proof}
a) Let $g\in\Z^n_\omega\backslash G_\fm$. By Lemma \ref{lem:gmtilde},
$A_gM_\fm=0$ and thus $\ph_\fm(A_g)M_\fm=0$.

b) Since $s_i\in G_\fm$, $\ph_\fm(a_{s_i}^*)b_i\in R_\fm\backslash\{0\}$
and by Lemma \ref{lem:NBabc}a) with $g=0$ and $h=s_i$ we have
$b_i\ph_\fm(a_{s_i}^*)\in R_\fm\backslash\{0\}$. So the $b_i$ are invertible.
%By Lemma \ref{lem:NBabc}b) and Lemma \ref{lem:dim},
%\[1\in R_\fm=\ph_\fm(A_0)=\ph_\fm(A_{s_i}A_{-s_i})=
%\ph_\fm(A_{s_i})\ph_\fm(A_{-s_i})=b_i R_\fm\ph_\fm(A_{-s_i})
%\]
%so $b_i$ has a right inverse and similarly it has a left inverse.
The relation (\ref{eq:Brel1}) follows from (\ref{eq:araMm1}).
Next we prove (\ref{eq:Brel2}).
From Lemma \ref{lem:NBabc}a) and Lemma \ref{lem:dim} we have
$\varphi(A_{s_i+s_j})=R_\fm b_ib_j$.
Switching $i$ and $j$ it follows that (\ref{eq:Brel2}) must hold for some
nonzero $\lambda_{ij}\in R_\fm$.

Finally we prove that (\ref{eq:bm1basis}) is a basis for $B_\fm^{(1)}$ over
$K$. Linear independence is clear. Let $g\in G_\fm$ and write
$g=\sum_i l_is_i$. By repeated use of Lemma \ref{lem:NBabc}b) we
obtain that
\[\ph_\fm(A_g)=\ph_\fm(A_{\sgn(l_1)s_1})^{|l_1|}\ldots
\ph_\fm(A_{\sgn(l_k)s_k})^{|l_k|}.\]
For $l_i=0$ the factor should be interpreted as $R_\fm$.
By Lemma \ref{lem:dim},
\[\ph_\fm(A_{\pm s_i})^l=R_\fm b_i^{\pm l}\]
for $l>0$ so using (\ref{eq:Brel1}) we get
\[\ph_\fm(A_g)=R_\fm b_1^{l_1}\ldots b_k^{l_k}.\]
The proof is finished.
\end{proof}

%
%
% A l g e b r a i c a l l y   c l o s e d   K
%
%
\subsection{Restricted case}\label{sec:restricted}
In this subsection we will assume that $K$ is algebraically closed.
Moreover we will assume that the $K$-algebra inclusion
$K\hookrightarrow R_\fm$ is
onto which is the case when $R$
is finitely generated as a $K$-algebra by the (weak) Nullstellensatz.
Then $\Z^n_\omega$ acts trivially on $R_\fm$.
The structure of $B_\fm^{(1)}$ given in Theorem \ref{thm:varphibmstruct1}
is then simplified in the following way.
\begin{cor}\label{cor:Bm1}
Let $k=\rank G_\fm$
and let $b_i=\varphi_\fm(a_{s_i})$ for $i=1,\ldots,k$ where
$\{s_1,\ldots,s_k\}$ is a $\Z$-basis for $G_\fm$.
Then $B_\fm^{(1)}$ is the $K$-algebra with
invertible generators $b_1,\ldots,b_k$ and the
relation
\[b_ib_j=\la_{ij}b_jb_i,\quad 1\le i,j\le k.\]
\end{cor}

Using the normal form of a skew-symmetric integral matrix
we will now show that $B_\fm^{(1)}$ can
be expressed as a tensor product of noncommutative tori.
Consider the matrix $(\la_{ij})_{1\le i,j\le k}$ from (\ref{eq:Brel2}).
\begin{claim}
If $B_\fm^{(1)}$ has a nontrivial irreducible
finite-dimensional representation, then all the $\la_{ij}$
are roots of unity.
\end{claim}
\begin{proof}
Indeed, let $N$ be a finite-dimensional simple module
over $B_\fm^{(1)}$ and let $i\in\{1,\ldots,k\}$.
Since $K$ is algebraically closed,
$b_i$ has an eigenvector $0\neq v\in N$ with eigenvalue $\mu$, say.
Since $b_i$ is invertible,
$\mu\neq 0$. Let $j\neq i$ and consider the vector $b_jv$. It is also nonzero,
since $b_j$ is invertible, and it is an eigenvector of $b_i$ with eigenvalue
$\lambda_{ij}\mu$. Repeating the process,
we obtain a sequence
\[\mu,\; \la_{ij} \mu,\; \la_{ij}^2 \mu,\;\ldots\]
of eigenvalues of $b_i$. Since $N$ is finite-dimensional, they cannot
all be pairwise distinct, and thus $\lambda_{ij}$ is a root of unity.
\end{proof}

For $\la\in K$, let $T_\la$ denote the $K$-algebra with two
invertible generators $a$ and $b$ satisfying $ab=\la ba$. $T_\la$
(or its $C^*$-analogue)
is sometimes referred to as a noncommutative torus.

\begin{thm} \label{thm:Bisotensor}
Let $k=\rank G_\fm$.
If all the $\la_{ij}$ in (\ref{eq:Brel2}) are roots of unity, then
there is a root of unity $\la$, an integer $r$ with
$0\le r\le\lfloor k/2\rfloor$
and positive integers $p_i, i=1,\ldots,r$ with $1=p_1|p_2|\ldots|p_r$
such that
\[B_\fm^{(1)}\simeq T_{\la^{p_1}}\otimes T_{\la^{p_2}}\otimes\cdots\otimes
T_{\la^{p_r}}\otimes L\]
where $L$ is a Laurent polynomial algebra over $K$ in $k-2r$ variables.
\end{thm}
\begin{proof}
If $k=1$, then $B_\fm^{(1)}\simeq K[b_1,b_1^{-1}]$ and $r=0$.
If $k>1$, let $p$ be the smallest positive integer such that
$\lambda_{ij}^p=1$ for all $i,j$. Using that $K$ is algebraically closed, we
fix a primitive $p$:th root of unity $\varepsilon\in K$.
Then there are integers $\theta_{ij}$ such that
\[\lambda_{ij}=\varepsilon^{\theta_{ij}}\]
and
\begin{equation}\label{eq:thetaskew}
\theta_{ji}=-\theta_{ij}.
\end{equation}
Equation (\ref{eq:thetaskew}) means that $\Theta=(\theta_{ij})$ is a $k\times k$
skew-symmetric integer matrix.
Next, consider a change of generators of the algebra $B_\fm^{(1)}$:
\begin{equation}\label{eq:change_genU}
b_i \mapsto b_i'=b_1^{u_{i1}}\cdots b_k^{u_{ik}}
\end{equation}
Such a change of generators can be done if we are given an invertible
$k\times k$ integer matrix $U=(u_{ij})$.
The new commutation relations are
\begin{align*}
b_i'b_j'&=b_1^{u_{i1}}\cdots b_k^{u_{ik}}b_1^{u_{j1}}\cdots b_k^{u_{jk}}=\\
&=\la_{11}^{u_{1i}u_{1j}}\ldots\la_{k1}^{u_{ki}u_{1j}}\cdot \ldots\\
&\quad\cdot \la_{1k}^{u_{1i}u_{kj}}\ldots\la_{kk}^{u_{ki}u_{kj}}\cdot b_j'b_i'=\\
&=\varepsilon^{\sum_{p,q}\theta_{pq}u_{pi}u_{qj}} b_j'b_i'
\end{align*}
Hence $\Theta'=U^T\Theta U$. By Theorem IV.1 in \cite{newman} there is a $U$
such that $\Theta'$ has the skew normal form
\[\begin{pmatrix}
0 & \theta_1   \\
-\theta_1 & 0  \\
&& 0 & \theta_2   \\
&&-\theta_2 & 0  \\
&&&& \ddots   \\
&&&&& 0 & \theta_r   \\
&&&&&-\theta_r & 0  \\
&&&&&&& \boldsymbol{0} \\
\end{pmatrix}\]
where $r\le\lfloor k/2 \rfloor$ is the rank of $\Theta$, the $\th_i$ are nonzero
integers, $\th_i|\th_{i+1}$ and $\boldsymbol{0}$ is a
$k-2r$ by $k-2r$ zero matrix.
Set $\la=\varepsilon^{\th_1}$ and $p_i=\th_i/\th_1$ for $i=1,\ldots,r$.
The claim follows.
\end{proof}

%
%
% S i m p l e   m o d u l e s   o v e r   n o n c o m m u t a t i v e   t o r i
%
%
%\subsection{Simple modules over noncommutative tori}
The following result, describing simple modules over the tensor
product of noncommutative tori, is more or less well-known, but we provide
a proof for convenience.
\begin{prop} \label{prop:tori1}
Let $M$ be a finite dimensional simple module over
\[T:=T_{\la_1}\otimes\cdots\otimes T_{\la_r},\]
where the $\la_i$ are roots of unity in $K$.
Then there are simple modules $M_i$ over $T_{\la_i}$ such that,
as $T$-modules,
\[M\simeq M_1\otimes\cdots\otimes M_r.\]
\end{prop}
\begin{proof}
Denote the generators of $T_{\la_i}$ by $a_i$ and $b_i$. We will
view $T_{\la_i}$ as subalgebras of $T$.
Since the elements $a_i, i=1,\ldots,r$ commute and
$M$ is finite dimensional and $K$ is algebraically closed,
there is a nonzero common eigenvector $w\in M$ of the $a_i$:
\begin{equation}\label{eq:thmtensor_ai}
a_i w = \mu_i w,\quad i=1,\ldots,r,
\end{equation}
where $\mu_i\in K^*$ because $a_i$ is invertible.
Let $n_i$ be the order of $\lambda_i$. Then $b_i^{n_i}$ acts as a scalar
by Schur's Lemma.
By simplicity of $M$, any element of $M$
has the form (using the commutation relations and (\ref{eq:thmtensor_ai}))
\begin{equation}\label{eq:thmtensor_elt}
\sum_{j\in\Z^r,\;0\le j_i<n_i} \rho_j b_1^{j_1}\ldots b_r^{j_r}\cdot w,
\end{equation}
where $\rho_j\in K$. This shows that
\[\dim_K M \le n_1\cdot\ldots\cdot n_r.\]
But the terms in (\ref{eq:thmtensor_elt}) all belong to different weight spaces
with respect to the commutative subalgebra generated by $a_1,\ldots,a_r$:
\[a_i\cdot b_1^{j_1}\ldots b_r^{j_r}w =
 \la_i^{j_i}\mu_i\cdot b_1^{j_1}\ldots b_r^{j_r} w,\quad i=1,\ldots,r,\]
and
\[(\la_1^{j_1}\mu_1,\ldots,\la_r^{j_r}\mu_r)\neq
(\la_1^{l_1}\mu_1,\ldots,\la_r^{l_r}\mu_r)\]
if $j,l\in\Z^r, 0\le j_i,l_i< n_i$ and $j\neq l$.
Hence by standard results they must be linearly independent. Thus
\begin{equation}\label{eq:thmtensor_dim}
\dim_K M = n_1\cdot\ldots\cdot n_r.
\end{equation}
Next, set $M_i = T_{\la_i}\cdot w$. Then
$M_i= \oplus_{j=0}^{n_i-1} Kb_i^j\cdot w$ and
\begin{equation}\label{eq:thmtensor_dim2}
\dim_K M_i = n_i.
\end{equation}
Finally, define
\[\psi: M_1\otimes\ldots\otimes M_r\to M\]
by
\[\psi(w\otimes\ldots\otimes w)=w\]
and by requiring that $\psi$ is a $T$-module homomorphism. This is possible
since $M_1\otimes\ldots\otimes M_r$ is generated by $w\otimes\ldots\otimes w$
as a $T$-module. Then
$\psi$ is surjective, since $M$ is simple. Also the dimensions
on both sides agree, so $\psi$ is an isomorphism of $T$-modules.
\end{proof}

%Since $T_\lambda$ is a GWA, Proposition \ref{prop:tori1} is a
%special case of the fact that
%any simple weight module over $A_1\otimes \ldots \otimes A_n$, where
%the $A_i$ are GWAs, is isomorphic to a tensor product
%$M_1\otimes\ldots\otimes M_n$ of simple weight modules $M_i$ over $A_i$.
%A proof of this can be found in \cite{?????}.

%
%
% D e s c r i p t i o n   o f   s i m p l e   m o d u l e s
%
%
\section{Explicit formulas for the induced modules}
\label{sec:explicit}
In this section we show explicitly how one can obtain
simple weight modules with no proper inner breaks
over a TGWA (equivalently over a TGWC by Proposition
\ref{prop:IM=0}) from the structure of its weight spaces
as $B(\omega)$-modules.

Since the $B(\omega)$-modules were described in the restricted
case in Subsection \ref{sec:restricted}, we obtain
in particular
a description of all simple weight modules over $A$
with no proper inner breaks and finite-dimensional
weight spaces if $R$ is finitely generated over
an algebraically closed field $K$.

\subsection{A basis for $M$}\label{sec:basis}
Let $\{v_i\}_{i\in I}$ be a basis for $M_\fm$ over $K$. By Lemma
\ref{lem:NBabc}a), $\tilde G_\fm$ is the union of some cosets in $\Z^n/G_\fm$. Let
$S\subseteq \Z^n$ be a set of representatives of these cosets.
For $g\in\tilde G_\fm$, choose $r_g\in R$ such that $a_g':=r_ga_g^*$ satisfies
$\ph_\fm(a_g')\ph_\fm(a_g)=1$.
\begin{thm}\label{thm:basis}
The set $C=\{a_gv_i \; |\; g\in S, i\in I\}$ is a basis for $M$ over $K$.
\end{thm}
\begin{proof}
First we show that $C$ is linearly independent over $K$.
Assume that
\[\sum_{g,i} \la_{gi} a_gv_i=0.\]
Then $\sum_i \la_{gi}a_gv_i=0$ for each
$g$ since the elements belong to different weight spaces. Hence
$0=a_g'\sum_i \la_{gi}a_gv_i= \sum_i \la_{gi}v_i$ for each $g$.
Since $v_i$ is a basis over $K$, all the $\la_{gi}$ must be zero.

Next we prove that $C$ spans $M$ over $K$. Since $M$ is simple and
$M_\fm\neq 0$,
\[M=AM_\fm=\sum_{g\in\Z^n}A_gM_\fm =\sum_{g\in\tilde G_\fm}A_gM_\fm=
\sum_{h\in S}\sum_{g\in h+G_\fm}A_gM_\fm =
\sum_{h\in S}A_hM_\fm\]
by Lemma \ref{lem:gmtilde} and Lemma \ref{lem:NBabc}c).
\end{proof}
\begin{cor}\label{cor:supp}
$\supp(M)=\{g(\fm)\; |\; g\in S\}$ and $g(\fm)\neq h(\fm)$ if $g,h\in S$, $g\neq h$.
\end{cor}
\begin{cor}\label{cor:dimM}
$\dim M=|S|\cdot\dim M_\fm$ with natural interpretation of $\infty$.
\end{cor}

%
% The action of A
%
\subsection{The action of $A$}
Our next step is to describe the action of the $X_i, Y_i$ on the basis $C$
for $M$.
Let $\zeta : \tilde G_\fm\to S$ be the function defined
by requiring $g-\zeta(g)\in G_\fm$.

\begin{thm}\label{thm:action}
Let $g\in S$ and let $v\in M_\fm$. Then
\[X_ia_gv=
\begin{cases}
a_h\cdot b_{g,i}v
& \text{if $g+e_i\in \tilde G_\fm$},\\
0& \text{otherwise},
\end{cases}\]
where $h=\zeta(g+e_i)$ and
\[b_{g,i}=(-h)(X_ia_g a_{g+e_i-h}'a_h')\cdot a_{g+e_i-h}\]
and
\[Y_ia_gv=
\begin{cases}
a_k\cdot c_{g,i}v
& \text{if $g-e_i\in \tilde G_\fm$},\\
0& \text{otherwise},
\end{cases}\]
where $k=\zeta(g-e_i)$
and
\[c_{g,i}=(-k)(Y_ia_g a_{g-e_i-k}'a_k')\cdot a_{g-e_i-k}.\]
\end{thm}
\begin{remark} Note that
\[\deg X_ia_g a_{g+e_i-h}'a_h'=\deg Y_ia_g a_{g-e_i-k}'a_k'=0\]
so the action of $\Z^n$ on these elements is well defined.
Thus we see that
$\deg b_{g,i}\in G_\fm$ and $\deg c_{g,i}\in G_\fm$, i.e. that
$b_{g,i}$ and $c_{g,i}$ belong to $B(\omega)$.
Therefore the action of these elements on a basis element
$v_i$ of $M_\fm$ can be determined if we know
the structure of $M_\fm$ as an $B(\omega)$-module. In the
restricted case this was described in
Subsection \ref{sec:restricted}.
Expanding the result in the basis $\{v_i\}$ again and
acting by $a_h$ or $a_k$ we obtain a linear
combination of basis elements from the set $C$.
\end{remark}
\begin{proof}
Assume $g+e_i\in\tilde G_\fm$. Let $h=\zeta(g+e_i)$. Then
\begin{align*}
X_ia_gv &= X_ia_g a_{g+e_i-h}'a_{g+e_i-h}v=\\
&=(X_ia_g a_{g+e_i-h}'a_h')a_h a_{g+e_i-h}v=\\
&=a_h\cdot (-h)(X_ia_g a_{g+e_i-h}'a_h')\cdot a_{g+e_i-h}v.
\end{align*}
If $g+e_i\notin\tilde G_\fm$, then $X_ia_gv=0$ by Lemma
\ref{lem:gmtilde}.

Assume $g-e_i\in\tilde G_\fm$. Let $k=\zeta(g-e_i)$. Then
\begin{align*}
Y_ia_gv &= Y_ia_g a_{g-e_i-k}'a_{g-e_i-k}v=\\
&=(Y_ia_g a_{g-e_i-k}'a_k')a_k a_{g-e_i-k}v=\\
&=a_k\cdot (-k)(Y_ia_g a_{g-e_i-k}'a_k')\cdot a_{g-e_i-k}v.
\end{align*}
If $g-e_i\notin\tilde G_\fm$, then $Y_ia_gv=0$ by Lemma
\ref{lem:gmtilde}.
\end{proof}

Note that we do not need the technical assumptions in the proof
of Theorem 1 in \cite{MT} under which the exact formulas for
simple weight modules were obtained.

%
%
% E x a m p l e :   q u a n t i z e d   W e y l   a l g e b r a
%
%
\section{Application to quantized Weyl algebras}
\label{sec:qweyl}
In this final part we will apply the methods developed in
the previous sections to the problem of describing representations
of the quantized Weyl algebra, defined in Section \ref{sec:examples}.
As mentioned there, it is naturally a TGWA.

First we find the isotropy group and the set $\tilde G_\fm$ expressed
as solution of systems of linear equations (see Proposition
\ref{prop:qweylsim2} and Proposition \ref{prop:tildeGm}). These sets
are directly related to the structure of the subalgebra $B(\omega)$
(Theorem \ref{thm:varphibmstruct1}) and the support of a module
(Corollary \ref{cor:supp}).

Then in Section \ref{sec:qweylrank2}
we give a complete classification
of all locally finite simple weight modules with no proper inner breaks
over a quantized Weyl algebra of rank two.
The parameters $q_1$ and $q_2$ are allowed
to be any numbers from $\C\backslash\{0,1\}$.
Example \ref{ex:qweylmodwithSPIB} shows that the assumption that
the modules have no proper inner breaks is not superfluous.

\subsection{The isotropy group and $\tilde G_\fm$}
Let $R=\C[t_1,\ldots,t_n]$ and fix
$\fm=(t_1-\al_1,\ldots,t_n-\al_n)\in\Max(R)$. Let $\omega$
be the orbit of $\fm$ under the action (\ref{eq:Znaction}) of $\Z^n$.
Set $[k]_q=\sum_{j=0}^{k-1}q^i$ for $k\in\Z$ and $q\in\C$.
Recall the definition (\ref{eq:qweylsigmadef}) of the
automorphisms $\si_i$ of $R$.
\begin{prop} \label{prop:qweylsim}
Let $(g_1,\ldots,g_n)\in\Z^n$. Then
\begin{align*}
\si_1^{g_1}\ldots\si_n^{g_n}(\fm)&=\\
&\Big([g_1]_{q_1}+q_1^{g_1}t_1-\al_1,\quad
[g_2]_{q_2}\big(1+(q_1-1)\al_1\big)+q_1^{g_1}q_2^{g_2}t_2-\al_2,\ldots\\
&\ldots,[g_j]_{q_j}\big(1+\sum_{r=1}^{j-1} (q_r-1)\al_r\big)+q_1^{g_1}\ldots q_j^{g_j}t_j-\al_j,\ldots\\
&\ldots,[g_n]_{q_n}\big(1+\sum_{r=1}^{n-1}(q_r-1)\al_r\big)+q_1^{g_1}\ldots q_n^{g_n}t_n-\al_n\Big).
\end{align*}
\end{prop}
\begin{proof} Induction.
\end{proof}

For notational brevity we set $\be_i=(q_i-1)\al_i$ and
$\ga_i=1+\be_1+\be_2+\ldots+\be_i$. We also set $\ga_0=1$. The numbers
$\ga_i$ will play an important role in the next statements.
By a \emph{$j$-break} we mean an ideal $\fn\in\Max(R)$ such that
$t_j\in \fn$.
\begin{cor} \label{cor:breaks}
For $j=1,\ldots,n$ we have
\[t_j\in \si_1^{g_1}\ldots\si_n^{g_n}(\fm)
\Longleftrightarrow \ga_j=q_j^{g_j}\ga_{j-1}.\]
Thus $\omega$ contains a $j$-break iff
$\ga_j=q_j^k\ga_{j-1}$ for some integer $k$.
\end{cor}
\begin{proof}
By Proposition \ref{prop:qweylsim},
\[t_j\in \si_1^{g_1}\ldots\si_n^{g_n}(\fm)\]
iff
\[[g_j]_{q_j}\big(1+\sum_{r=1}^{j-1}
 (q_r-1)\al_r\big)=\al_j.\]
Multiply both sides with $q_j-1$ to get
\[(q_j^{g_j}-1)(1+\be_1+\ldots+\be_{j-1})=\be_j.\]
\end{proof}

The next Proposition describes the isotropy subgroup $\Z^n_\omega$
defined in (\ref{eq:Znga}).
\begin{prop} \label{prop:qweylsim2}
We have
\begin{equation} \label{eq:qweylsim2pf}
\Z^n_\omega=\{g\in\Z^n\; |\; (q_1^{g_1}\ldots q_j^{g_j}-1)\ga_j=0\;\forall j=1,\ldots,n\}.
\end{equation}
\end{prop}
\begin{proof} From Proposition \ref{prop:qweylsim},
$\si_1^{g_1}\ldots\si_n^{g_n}(\fm)=\fm$ iff
\begin{align*}
\al_1&=[g_1]_{q_1}+q_1^{g_1}\al_1\\
\al_2&=[g_2]_{q_2}\big(1+(q_1-1)\al_1\big)+q_1^{g_1}q_2^{g_2}\al_2\\
&\vdots\\
\al_n&=[g_n]_{q_n}\big(1+(q_1-1)\al_1+\ldots+(q_{n-1}-1)\al_{n-1}\big)+
q_1^{g_1}\ldots q_n^{g_n}\al_n
\end{align*}
Multiply the $i$:th equation by $q_i-1$. Then the system can be written
\begin{align*}
\be_1&=q_1^{g_1}-1+q_1^{g_1}\be_1\\
\be_2&=(q_2^{g_2}-1)(1+\be_1)+q_1^{g_1}q_2^{g_2}\be_2\\
&\vdots\\
\be_n&=(q_n^{g_n}-1)(1+\be_1+\ldots+\be_{n-1})+q_1^{g_1}\ldots q_n^{g_n}\be_n
\end{align*}
or equivalently
\begin{align*}
1+\be_1&=q_1^{g_1}(1+\be_1)\\
1+\be_1+\be_2&=q_2^{g_2}(1+\be_1)+q_1^{g_1}q_2^{g_2}\be_2\\
&\vdots\\
1+\be_1+\ldots+\be_n&=q_n^{g_n}(1+\be_1+\ldots+\be_{n-1})+q_1^{g_1}\ldots q_n^{g_n}\be_n
\end{align*}
Now for $i$ from $1$ to $n-1$, replace the expression $1+\be_1+\ldots+\be_i$ in the
right hand side of the $i+1$:th equation by the right hand side
of the $i$:th equation. After simplification,
the claim follows.
\end{proof}
Note that it follows from (\ref{eq:qweylsim2pf})
that the subgroup
\begin{equation}\label{eq:Q}
Q=\{g\in\Z^n\; |\; q_j^{g_j}=1\text{ for }j=1,\ldots,n\}
\end{equation}
of $\Z^n$ is always contained in $\Z^n_\omega$ for any orbit
$\omega$. Moreover $\Z^n_\omega=Q$ if $\omega$ (viewed as
a subset of $\C^n$) does not intersect the union of the
hyperplanes in $\C^n$ defined by the equations
$1+(q_1-1)x_1+\ldots+(q_j-1)x_j=0$ ($1\le j\le n$).
%Of course the group $Q$ can be trivial. This is the
%case for example when all the $q_j$ are positive real
%numbers.

Another case of interest is when for any $j$,
$q_1^{g_1}\ldots q_j^{g_j}=1$ implies $g_1=\ldots =g_j=0$.
If for instance the $q_j$ are pairwise distinct prime numbers this hold.
Then $\Z^n_\omega=\{0\}$ unless $1+\be_1+\ldots+\be_j=0$
for all $j$, i.e. unless $\omega$ contains the point
\[\fn_0=(t_1-(1-q_1)^{-1},t_2,\ldots,t_n).\]
So in this very special case we have $\omega=\{\fn_0\}$ and
$\Z^n_\omega=\Z^n$.

%---------------------------------------------------------------------------

We now turn to the set $\tilde G_\fm$ defined in (\ref{eq:tildeGmandGm})
which can here be described explicitly in terms of $\fm$ in the following way.
\begin{prop} \label{prop:tildeGm}
\[\tilde G_\fm=\tilde G_\fm^{(1)}\times \ldots\times
\tilde G_\fm^{(n)},\]
where
\begin{align*}
\tilde G_\fm^{(j)}&=\{k\ge 0\; |\; \ga_j\neq q_j^i\ga_{j-1}
\;\forall i=0,1,\ldots,k-1\}\cup\\
&\quad\cup\{k<0\; |\; \ga_j\neq q_j^i\ga_{j-1}\;\forall i=-1,-2,\ldots,k\}.
\end{align*}
\end{prop}
\begin{proof} From the relations of the algebra
follows that the subspace spanned by the words in $A_g$
is one-dimensional. Thus $g\in \tilde G_\fm$ iff
\begin{equation}\label{eq:qweyltildeGmpf}
Z_n^{-g_n}\ldots Z_1^{-g_1}Z_1^{g_1}\ldots Z_n^{g_n}\notin\fm
\end{equation}
where $Z_i^k=X_i^k$ if $k\ge 0$ and $Z_i^k=Y_i^{-k}$ if $k<0$.
Since $\si_i(t_j)=t_j$ for $j<i$, (\ref{eq:qweyltildeGmpf}) is equivalent to
\[Z_n^{-g_n}Z_n^{g_n}\ldots Z_1^{-g_1}Z_1^{g_1}\notin\fm.\]
Since $\fm$ is prime, this holds iff $Z_j^{-g_j}Z_j^{g_j}\notin\fm$
for each $j$. If $g_j=0$ this is true. If $g_j>0$ we have
\[Z_j^{-g_j}Z_j^{g_j}=Y_j^{g_j}X_j^{g_j}=Y_j^{g_j-1}X_j^{g_j-1}\si_j^{-g_j+1}(t_j)=\ldots
=t_j\si_j^{-1}(t_j)\ldots\si_j^{-g_j+1}(t_j),\]
while if $g_j<0$
\[Z_j^{-g_j}Z_j^{g_j}=X_j^{-g_j}Y_j^{-g_j}=X_j^{-g_j-1}Y_j^{-g_j-1}\si_j^{-g_j}(t_j)=\ldots
=\si_j(t_j)\ldots\si_j^{-g_j}(t_j).\]
Since $\fm$ is prime,
$g\in \tilde G_\fm$ iff for all $j=1,\ldots,n$
\[
t_j\notin\si_j^{i}(\fm),\quad i=0,\ldots,g_j-1 \text{ if }
g_j\ge 0,
\]
and
\[
t_j\notin\si_j^{i}(\fm),\quad i=-1,-2\ldots,g_j \text{ if }
g_j< 0.
\]
The claim now follows from Corollary \ref{cor:breaks}.
\end{proof}

%------------------------------------------------------------

\begin{cor} If $\{1,\al_1,\al_2,\ldots,\al_n\}$ is linearly
independent over $\Q(q_1,\ldots,q_n)$, then $\tilde G_\fm=\Z^n$.
\end{cor}

%
%
%  D e g e n e r a t e   o r b i t s   i n   r a n k   t w o
%
%
\subsection{Description of simple weight modules over rank two algebras}
\label{sec:qweylrank2}
Assume from now on that $A$ is a quantized Weyl algebra of rank two.
In this section we will obtain a list of all locally finite simple
weight $A$-modules with no proper inner breaks.

We consider first some families of ideals in $\Max(R)$.
Define for $\la\in\C$,
\begin{align*}
\fn_\la^{(1)} &= \big(t_1-(1-\la)(1-q_1)^{-1},t_2-\la(1-q_2)^{-1}\big),\\
\fn_\la^{(2)} &= \big(t_1-(1-q_1)^{-1},t_2-\la\big),
\end{align*}
and set $\fn_0=\fn_0^{(1)}=\fn_0^{(2)}$.
The following lemma will be useful.
\begin{lem}\label{lem:si1ksi2l}
For $\la\in\C$ and integers $k,l$ we have
\begin{align}
\si_1^k\si_2^l(\fn_\la^{(1)})&=\fn_{\la q_1^{-k}}^{(1)},\label{eq:qwsl1}\\
\si_1^k\si_2^l(\fn_\la^{(2)})&=\fn_{\la q_1^{-k}q_2^{-l}}^{(2)}.
\label{eq:qwsl2}
\end{align}
\end{lem}
\begin{proof} Follows from Proposition \ref{prop:qweylsim} or by direct
calculation using the definition (\ref{eq:qweylsigmadef}) of the $\si_i$.
\end{proof}

The following example shows the existence of locally finite simple weight
modules $M$ over $A$ which have some proper inner breaks.
% example with proper inner breaks
\begin{example} \label{ex:qweylmodwithSPIB}
%Let $A$ be a quantized Weyl algebra of rank two and
Assume that $q_1\la_{12}$ is a root of unity of order $r$. Let $M$ be
a vector space of dimension $r$ and let $\{v_0, v_1,\ldots, v_{r-1}\}$
be a basis for $M$. Define an action of $A$ on $M$ as follows.
\begin{align*}
X_1v_k&=\begin{cases}v_{k+1},&k<r-1\\v_0,&k=r-1\end{cases} & X_2v_k&=(q_1\la_{12})^{-k}v_k\\
Y_1v_k&=\begin{cases}(1-q_1)^{-1}v_{k-1},&k>0\\(1-q_1)^{-1}v_{r-1},&k=0\end{cases} & Y_2v_k&=0
\end{align*}
It is easy to check that (\ref{eq:qweylrel1})--(\ref{eq:qweylrel3}) hold so this
defines a module over $A$. It is immediate that $M=M_\fm$
where $\fm=\fn_0=(t_1-(1-q_1)^{-1},t_2)$ so $M$ is a weight module and $M$ is simple
by standard arguments. However, recalling Definition \ref{dfn:NPIB},
$M$ has some proper inner breaks in the sense
that $\fm\in\supp(M)$, $X_2M_\fm\neq 0$ but $Y_2X_2M_\fm=0$.
\end{example}

We will describe the isotropy groups of the different ideals in $\Max(R)$.
Let $K_1$ and $K_2$ denote the kernels of the group homomorphisms from $\Z\times \Z$ to
the multiplicative group $\C\backslash\{0\}$ which map $(k,l)$ to $q_1^k$ and $q_1^kq_2^l$
respectively. Then $Q=K_1\cap K_2$ where $Q$ was defined in (\ref{eq:Q}).
For $\fm\in\Max(R)$, recall that $\Z^2_\fm=\{g\in\Z^2\; |\; g(\fm)=\fm\}$.
The following corollary describes the isotropy group
$\Z^2_\fm$ of any $\fm\in\Max(R)$.
\begin{cor}\label{cor:Z2fm}
Let $\la\in\C\backslash\{0\}$ and $\fn\in\Max(R)\backslash
\{\fn_{\mu}^{(i)}\; |\; \mu\in\C, i=1,2\}$. Then we have the following
equalities in the lattice of subgroups of $\Z^2$.
\[\xymatrix{
  & \Z_{\fn_0}^2=\Z^2  \ar@{-}[dl] \ar@{-}[dr] \\
\Z_{\fn_\la^{(1)}}^2=K_1 \ar@{-}[dr]  && \Z_{\fn_\la^{(2)}}^2=K_2 \ar@{-}[dl] \\
& \Z_{\fn}^2=Q        }\]
\end{cor}
\begin{proof}
The family of ideals $\{\fn_{\la}^{(1)} \; |\; \la\in\C\}$ are precisely those
for which $\ga_2=0$. And $\{\fn_{\la}^{(2)} \; |\; \la\in\C\}$ are exactly
those such that $\ga_1=0$. Thus the claim follows from
Proposition \ref{prop:qweylsim2}.
\end{proof}

Let $M$ be a simple weight $A$-module with no proper inner breaks
and finite dimensional weight spaces,
$\fm=(t_1-\al_1,t_2-\al_2)\in\supp M$ and let $\omega$ be the orbit of $\fm$.
We consider four main cases separately:
$\fm=\fn_0$, $\fm=\fn^{(1)}_\la$ for some $\la\neq 0$,
$\fm=\fn^{(2)}_\la$ for some $\la\neq 0$ and
$\fm\notin\{\fn_{\mu}^{(i)}\; |\; \mu\in\C, i=1,2\}$.
Some of these cases will contain subcases.
In each case we will proceed along the following steps, which also illustrate
the procedure for a general TGWA.
\begin{enumerate}
\item Find the sets $\Z^n_\fm$ and $\tilde G_\fm$
using Corollary \ref{cor:Z2fm} and Proposition \ref{prop:tildeGm}.
Write down $G_\fm=\Z^n_\fm\cap\tilde G_\fm$ and
choose a basis $\{s_1,\ldots,s_k\}$ for $G_\fm$ over $\Z$.
\item For each $g\in\tilde G_\fm$, choose a word $a_g$ of degree $g$ such that
$a_g^*a_g\notin\fm$.
\item Using Corollary \ref{cor:Bm1}, describe $B_\fm^{(1)}$ and the
 finite-dimensional simple $B_\fm^{(1)}$-module $M_\fm$.
\item Choose a set of representatives $S$ for $\tilde G_\fm/G_\fm$.
 By Theorem \ref{thm:basis} we know then a basis $C$ for $M$.
\item Calculate the action of $X_i$, $Y_i$ on the basis using
either relations (\ref{eq:qweylrel1})--(\ref{eq:qweylrel3}) or Theorem \ref{thm:action}.
\end{enumerate}

We will use the following notation: $Z_j^k=X_j^k$ if $k\ge 0$ and
$Z_j^k=Y_j^{-k}$ if $k<0$.
Note that the $k$ in $Z_j^k$ should only be regarded as an upper index, not
as a power. The choice of $a_g$ in step two above is more or less
irrelevant for a quantized Weyl algebra because each $A_g$ is
one-dimensional. Therefore we will always choose
$a_g=Z_1^{g_1}Z_2^{g_2}$ where $g=(g_1,g_2)$.

%
%
%    F i r s t   c a s e
%
%
\subsection{The case $\fm=\fn_0$}
Here $\al_1=(1-q_1)^{-1}$, $\al_2=0$ so that $\ga_1=\ga_2=0$.
By Corollary \ref{cor:Z2fm} we have $\Z_\fm^2=\Z^2$ and
from Proposition \ref{prop:tildeGm} one obtains
that $\tilde G_\fm=\Z\times\{0\}$. Thus
$G_\fm=\Z\times\{0\}=\Z\cdot s_1$ with $s_1=(1,0)$. Since $G_\fm$ has
rank one, Corollary \ref{cor:Bm1} implies that $B_\fm^{(1)}$ is isomorphic to
the Laurent polynomial algebra $\C[T,T^{-1}]$ in one variable.
Therefore $M_\fm$ is one-dimensional, say $M_\fm=\C v_0$ and
$b_1=\varphi_\fm(Z_1^1)=\varphi_\fm(X_1)$, hence $X_1$, acts in
$M_\fm$ as some nonzero scalar $\rho$. And
\[Y_1v_0=\rho^{-1}Y_1X_1v_0=\rho^{-1}(1-q_1)^{-1}v_0.\]
Here $S=\{(0,0)\}$ and $C=\{v_0\}$ is a basis for $M$ with the following
action:
\begin{align}
X_1v_0&=\rho v_0, & X_2v_0&=0,\label{eq:cl_n0}\\
Y_1v_0&=\rho^{-1}(1-q_1)^{-1}v_0, &Y_2v_0&=0.\nonumber
\end{align}
That $Z_2^{\pm 1}v_0=0$ follows from Theorem \ref{thm:action}
since $(0,\pm 1)\notin \tilde G_\fm$.

%
%
%   S e c o n d   c a s e
%
%
\subsection{The case $\fm=\fn_\la^{(1)}$, $\la\neq 0$}
Here $\al_1=(1-\la)(1-q_1)^{-1}$ and $\al_2=\la(1-q_1)^{-1}$ so
$\ga_1=\la$ and $\ga_2=0$. By Proposition \ref{prop:tildeGm},
$\tilde G_\fm^{(2)}=\Z$ and
\[
\tilde G_\fm^{(1)}=\{k\ge 0\; |\; \la\neq q_1^i
\;\forall i=0,1,\ldots,k-1\}
\cup\{k<0\; |\; \la\neq q_1^i\;\forall i=-1,-2,\ldots,k\}.
\]

We consider four subcases according to whether $\omega$
contains a $1$-break or not and whether $q_1$ is a root of unity
or not.
\subsubsection{The case $\fm=\fn_\la^{(1)}$, $\la\neq 0$, $\omega$
contains a $1$-break and $q_1$ is a root of unity}
By Corollary \ref{cor:breaks} $\la=q_1^k$ for some $k\in\Z$.
Let $o_1$ be the order of $q_1$. Then $\Z_\fm^2=K_1=(o_1\Z)\times\Z$.
We can further assume that $k\in\{0,1,\ldots,o_1-1\}$.

Note that $X_1^kM_\fm\neq 0$ because $\deg X_1^k=(k,0)\in\tilde G_\fm$
so $Y_1^kX_1^k\notin\fm$. Hence $\si_1^k(\fm)\in\supp(M)$. By
Lemma \ref{lem:si1ksi2l}, $\si_1^k(\fm)=\fn_{q_1^kq_1^{-k}}^{(1)}=
\fn_1^{(1)}$. We can thus change notation and let $\fm=\fn_1^{(1)}$.
Then by Proposition \ref{prop:tildeGm} we have
\[\tilde G_\fm=\{0,-1,-2,\ldots,-o_1+1\}\times\Z.\]
And $G_\fm=\tilde G_\fm\cap \Z^2_\fm=\{0\}\times\Z$.
By Corollary \ref{cor:Bm1}, $B_\fm^{(1)}$ is a
Laurent polynomial algebra in one variable. Thus $M_\fm$ is
one dimensional with a basis vector, say $v_0$. $X_2$
acts by some nonzero scalar $\rho$ on $v_0$ and $Y_2X_2v_0=(1-q_2)^{-1}v_0$.
$X_1$ and $Y_1^{o_1}$ act as zero on $M_\fm$ by Lemma \ref{lem:gmtilde}
because their degrees $(1,0)$ and $(-o_1,0)$ does not belong to $\tilde G_\fm$.

As a set of representatives for $\tilde G_\fm / G_\fm$ we choose
\[S=\{(0,0),(-1,0),(-2,0),\ldots,(-o_1+1,0)\}.\]
By Corollary \ref{cor:supp} we obtain that
\[\supp(M)=\{\fn_1^{(1)},\fn_{q_1^{-1}}^{(1)},
\ldots,\fn_{q_1^{-o_1+1}}^{(1)}\}.\]
By \ref{thm:basis}, the set
\[C=\{v_j:=Y_1^jv_0\; |\; j=0,1,\ldots,o_1-1\}\]
is a basis for $M$. The following picture shows the support of the module
and how the $X_i$ act on it. Since the $Y_i$ just act in the opposite
direction of the $X_i$ we do not draw their arrows.
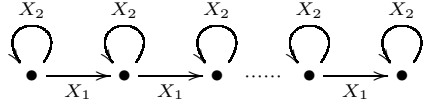
\[\xymatrix{
\bullet \ar[r]_{X_1} \ar@(ur,ul)_{X_2}&
\bullet \ar[r]_{X_1}  \ar@(ur,ul)_{X_2}&
\bullet \ar@{}[r]|{\cdots\cdots}  \ar@(ur,ul)_{X_2}&
\bullet \ar[r]_{X_1} \ar@(ur,ul)_{X_2}&
\bullet \ar@(ur,ul)_{X_2}
}\]
Using Lemma \ref{lem:si1ksi2l},
\[X_1v_j=X_1Y_1^jv_0=Y_1^{j-1}\si_1^j(t_1)v_0=[j]_{q_1}v_{j-1}\]
and from relations (\ref{eq:qweylrel1})--(\ref{eq:qweylrel3}) follow that
\[X_2v_j=q_1^j\la_{12}^jY_1^jX_2v_0=\rho\la_{12}^jq_1^jv_j,\]
\[Y_2v_j=\la_{21}^jY^jY_2v_0=(1-q_2)^{-1}\rho^{-1}\la_{21}^jv_j.\]
Thus the action on the basis $\{v_0,\ldots,v_{o_1-1}\}$ is
\begin{equation}\label{eq:cl_n1i}\begin{split}
X_1v_j&=
\begin{cases}
0,&j=0,\\
[j]_{q_1}v_{j-1},&0<j\le o_1-1,
\end{cases}\\
Y_1v_j&=
\begin{cases}
v_{j+1},&0\le j<o_1-1,\\
0,&j=o_1-1,
\end{cases}\\
X_2v_j&=\rho\la_{12}^jq_1^jv_j,\\
Y_2v_j&=(1-q_2)^{-1}\rho^{-1}\la_{21}^jv_j.
\end{split}\end{equation}

\subsubsection{The case $\fm=\fn_\la^{(1)}$, $\la\neq 0$, $\omega$
contains a $1$-break and $q_1$ is not a root of unity}
Now there is a unique integer $k\in\Z$ such that $\la=q_1^k$.
If $k\ge 0$, then $\tilde G_\fm^{(1)}$ is the
set of all integers $\le k$ while if $k<0$,
then $\tilde G_\fm^{(1)}$ is all integers
$\ge k+1$.

If $k\ge 0$,
$X_1^kM_\fm\neq 0$ because $(k,0)\in\tilde G_\fm$
so $Y_1^kX_1^k\notin\fm$.
Therefore $\si_1^k(\fm)=\fn_1^{(1)}\in\supp(M)$.
We change notation and let $\fm=\fn_1^{(1)}$.
Then $\tilde G_\fm^{(1)}=\{\ldots,-2,-1,0\}$
and $G_\fm=\{0\}\times\Z$.
We choose $S=\{(i,0)\; |\; i\le 0\}$.
$Y_2X_2=(1-q_2)^{-1}$ on $M_\fm$
so $M_\fm=\C v_0$, for a basis vector $v_0$, and $X_2v_0=\rho v_0$
for some $\rho\in\C^*$.
The set
$C=\{v_j:=Y_1^jv_0\; |\; j\le 0\}$
is a basis for $M$ and we have the following picture of $\supp(M)$.
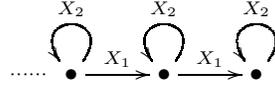
\[\xymatrix{
 \ar@{}[r]|{\cdots \cdots }&
\bullet \ar[r]^{X_1} \ar@(ur,ul)_{X_2}&
\bullet \ar[r]^{X_1} \ar@(ur,ul)_{X_2}&
\bullet \ar@(ur,ul)_{X_2}
}\]
One easily obtains the following action on the basis
$\{v_j\;|\; j\le 0\}$:
\begin{equation}\label{eq:cl_n1iihi}\begin{split}
X_1v_j&=
\begin{cases}
0,&j=0,\\
[j]_{q_1}v_{j-1},&j\ge 1,
\end{cases}\\
Y_1v_j&=v_{j+1},\\
X_2v_j&=\rho\la_{12}^jq_1^jv_j,\\
Y_2v_j&=(1-q_2)^{-1}\rho^{-1}\la_{21}^jv_j.
\end{split}\end{equation}

The case $k<0$ is analogous and yields a lowest weight
representation with $\fm=\fn_{q_1^{-1}}^{(1)}$ as its
lowest weight. A basis for $M$ is then
\[C=\{v_j:=X_1^jv_0\; |\; j\ge 0\},\]
where $M_\fm=\C v_0$ and the action is given by
\begin{equation}\label{eq:cl_n1iilo}\begin{split}
X_1v_j&=v_{j+1},\\
Y_1v_j&=\begin{cases}
0,&j=0,\\
[-j]_{q_1}v_{j-1},&j>0,
\end{cases}\\
X_2v_j&=(q_1\la_{12})^{-j}\rho v_j,\\
Y_2v_j&=\la_{12}^j(1-q_2)^{-1}\rho^{-1}v_j.
\end{split}\end{equation}

\subsubsection{The case $\fm=\fn_\la^{(1)}$, $\la\neq 0$, $\omega$
contains no $1$-break and $q_1$ is a root of unity}
By Corollary \ref{cor:breaks}, $\la\neq q_1^k$
for all $k\in\Z$. So by Proposition \ref{prop:tildeGm},
$\tilde G_\fm=\Z^2$. $G_\fm=(o_1\Z)\times\Z$
and we can choose $S=\{0,1,\ldots,o_1-1\}\times\{0\}$.
From
\[X_1^{o_1}X_2=(q_1\la_{12})^{o_1}X_2X_1^{o_1}=
\la_{12}^{o_1}X_2X_1^{o_1}\]
and Corollary \ref{cor:Bm1} follows that
$B_\fm^{(1)}\simeq T_{\la_{12}^{o_1}}$.
It can only have finite-dimensional irreducible
representations if $\la_{12}^{o_1}$ is a root of unity.
Assuming this, any such representation is $r$-dimensional,
where $r$ is the order of $\la_{12}^{o_1}$, and is
parametrized by $\C^*\times \C^*\ni (\rho,\mu)$
with basis
\[M_\fm=\spn\{v_j:=X_2^jv_0\; |\; j=0,1,\ldots,r-1\},\]
where $X_1^{o_1}v_0=\rho v_0$ and relations
\begin{align*}
X_1^{o_1}v_j&=\la_{12}^{o_1j}\rho v_j,\\
X_2v_j&=
\begin{cases}
v_{j+1},&0\le j<r-1,\\
\mu v_0, &j=r-1.
\end{cases}
\end{align*}
Therefore by Theorem \ref{thm:basis},
\[M=\spn\{w_{ij}=X_1^iv_j\; |\; 0\le i<o_1, 0\le j<r\}.\]
Using the commutation relations and the formulas in Lemma
\ref{lem:si1ksi2l} we can write down the action as
follows.

\begin{equation}\label{eq:cl_n1iii}\begin{split}
X_1w_{ij}&=
\begin{cases}
w_{i+1,j},&0\le i<o_1-1,\\
\la_{12}^{o_1j}\rho w_{0,j},& i=o_1-1,
\end{cases}\\
Y_1w_{ij}&=
\begin{cases}
(1-\la)(1-q_1)^{-1}\la_{12}^{-o_1j}\rho^{-1}w_{o_1-1,j},&i=0,\\
(1-\la q_1^{-i})(1-q_1)^{-1}w_{i-1,j},& 0<i\le o_1-1,
\end{cases}\\
X_2w_{ij}&=
\begin{cases}
q_1^{-i}\la_{21}^iw_{i,j+1},& 0\le j<r-1,\\
q_1^{-i}\la_{21}^i\mu w_{i,0},& j=r-1,
\end{cases}\\
Y_2w_{ij}&=
\begin{cases}
\la_{12}^i\mu^{-1}\la(1-q_2)^{-1}w_{i,r-1},&j=0,\\
\la_{12}^i\la(1-q_2)^{-1}w_{i,j-1},&0<j\le r-1.
\end{cases}
\end{split}\end{equation}
The action can be illustrated in the following way.
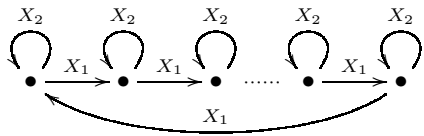
\[\xymatrix{
\bullet \ar[r]^{X_1} \ar@(ur,ul)_{X_2}&
\bullet \ar[r]^{X_1}  \ar@(ur,ul)_{X_2}&
\bullet \ar@{}[r]|{\cdots\cdots}  \ar@(ur,ul)_{X_2}&
\bullet \ar[r]^{X_1} \ar@(ur,ul)_{X_2}&
\bullet \ar@(ur,ul)_{X_2} \ar@(dl,dr)[llll]_{X_1}
}\]

\subsubsection{The case $\fm=\fn_\la^{(1)}$, $\la\neq 0$, $\omega$
contains no $1$-break and $q_1$ is not a root of unity}
By Corollary \ref{cor:breaks}, $\la\neq q_1^k$
for all $k\in\Z$. Now $\Z^2_\fm=\{0\}\times\Z$ so
$G_\fm=\{0\}\times\Z$. $M_\fm$ is one-dimensional with basis
$v_0$, say, and $X_2=\rho$ on $M_\fm$ while $Y_2X_2=\la(1-q_2)^{-1}\neq 0$
on $M_\fm$. We choose $S=\Z\times\{0\}$. Then a basis for $M$ is
\[C=\{v_j:=X_1^jv_0\; |\; j\ge 0\}\cup
 \{v_j:=\zeta_j Y_1^{-j}v_0\; |\; j<0\},\]
where we determine $\zeta_j$ by requiring that $X_1v_j=v_{j+1}$
for all $j$. Explicitly we have for $j<0$,
\[\zeta_j=\frac{(1-q_1)^{-j}}{(1-\la q_1^{-j})(1-\la q_1^{-j-1})\ldots
(1-\la q_1)}.\]
Using the commutation relations and the formulas in Lemma
\ref{lem:si1ksi2l} we get the action on $M=\spn\{v_j\; |\; j\in\Z\}$.
\begin{align}\label{eq:cl_n1iv}
X_1v_j&=v_{j+1}, &
X_2v_j&=q_1^{-j}\la_{12}^{-j}\rho v_j,\\
Y_1v_j&=\frac{1-\la q_1^{-j+1}}{1-q_1}v_{j-1},&
Y_2v_j&=\la_{12}^j\la(1-q_2)^{-1}\rho^{-1}v_j, \nonumber
\end{align}
and a corresponding diagram
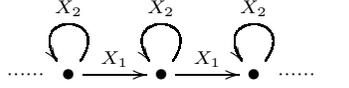
\[\xymatrix{
 \ar@{}[r]|{\cdots\cdots} &
\bullet \ar[r]^{X_1}  \ar@(ur,ul)_{X_2}&
\bullet \ar[r]^{X_1}  \ar@(ur,ul)_{X_2}&
\bullet \ar@{}[r]|{\cdots\cdots}  \ar@(ur,ul)_{X_2}&
}.\]

%
%
%   T h i r d   c a s e
%
%
\subsection{The case $\fm=\fn_\la^{(2)}$, $\la\neq 0$}
Here $\ga_1=0$ while $\ga_2=\la(q_2-1)$.
By Corollary \ref{cor:breaks}, $\omega$ does not contain
any breaks. We have $\tilde G_\fm=\Z^2$ and $G_\fm=\Z^2_\fm=K_2$.

We will need some lemmas in order to proceed.
\begin{lem} \label{lem:Z1kZ2lswap} For $k,l\in\Z$ we have
\begin{equation}\label{eq:Z1kZ2lswap}
Z_1^kZ_2^l=q_1^{k\bar l}\la_{12}^{kl}Z_2^lZ_1^k,
\end{equation}
where $\bar l=\max\{0,l\}$.
\end{lem}
\begin{proof}
Relations (\ref{eq:qweylrel1})--(\ref{eq:qweylrel3}) can be rewritten
in the more compact form
\[Z_1^kZ_2^l=q_1^{k\delta_{l,1}}\la_{12}^{kl}Z_2^lZ_1^k,\quad
k,l=\pm 1,\]
where $\delta_{l,1}$ is the Kronecker symbol.
After repeated application of this, (\ref{eq:Z1kZ2lswap}) follows.
\end{proof}

By Lemma \ref{lem:si1ksi2l} we have for $k,l\in\Z$,
\begin{align}
\si_1^k\si_2^l(t_1)&=(1-q_1)^{-1}\mod \fm,
\label{eq:si1ksi2ligen}
\\
\si_1^k\si_2^l(t_2)&=\la q_1^kq_2^l\mod \fm.
\label{eq:si1ksi2ligen2}
\end{align}

\begin{lem}\label{lem:ZikZiladd}
Let $k,l\in\Z$ and let $m=\min\{|k|,|l|\}$. Then,
as operators on $M_\fm$, we have
\begin{align}\label{eq:Z1kZ1ladd}
Z_1^kZ_1^l&=\begin{cases}
Z_1^{k+l},& kl\ge 0,\\
(1-q_1)^{-m}Z_1^{k+l},&kl<0,
\end{cases}\\
\label{eq:Z2kZ2ladd}
Z_2^kZ_2^l&=\begin{cases}
Z_2^{k+l},& kl\ge 0,\\
\la^m q_2^{(1-2l+(\sgn l)m)m/2}Z_2^{k+l},&kl<0.
\end{cases}
\end{align}
\end{lem}
\begin{proof}
Direct calculation using (\ref{eq:si1ksi2ligen})
and (\ref{eq:si1ksi2ligen2}). For example if $k>0$ and $l<0$ we have
\begin{align*}
Z_2^kZ_2^l&=X_2^kY_2^{-l}=X_2^{k-1}\si_2(t_2)Y_2^{-l-1}=\\
&=X_2^{k-1}Y_2^{-l-1}\si_2^{-l}(t_2)=
X_2^{k-1}Y_2^{-l-1}\la q_2^{-l}=\ldots=\\
&=\la q_2^{-l}\la q_2^{-l-1}\ldots \la q_2^{-l-(m-1)}Z_2^{k+l}=\\
&=\la^m q_2^{-lm-m(m-1)/2}Z_2^{k+l}.
\end{align*}
\end{proof}

\begin{lem}\label{lem:ZikZilswap}
Let $k,l\in\Z$ and let $m=\min\{|k|,|l|\}$. Then,
as operators on $M_\fm$,
\begin{equation}\label{eq:Z1kZ1lswap}
Z_1^kZ_1^l=Z_1^lZ_1^k,
\end{equation}
and
\begin{equation}\label{eq:Z2kZ2lswap}
Z_2^kZ_2^l=c(k,l)Z_2^lZ_2^k,
\end{equation}
where
\begin{equation}\label{eq:ckl}
c(k,l)=\begin{cases}
1,& kl\ge 0,\\
q_2^{(k-l)m-(\sgn k-\sgn l)m^2/2},&kl<0.
\end{cases}
\end{equation}
\end{lem}
\begin{proof} Follows directly from Lemma \ref{lem:ZikZiladd}.
\end{proof}

\begin{lem} \label{lem:rg}
Let $g=(g_1,g_2)\in\Z^2=\tilde G_\fm$ and set $r_g=\varphi_\fm(a_g^*a_g)^{-1}$
where $\varphi_\fm$ is the projection $R\to R/\fm=K$. Then
\begin{equation}\label{eq:rg}
r_g=(1-q_1)^{|g_1|}(\la^{-1} q_2^{(g_2-1)/2})^{|g_2|}
\end{equation}
and $(a_g)^{-1}=r_ga_g^*=r_gZ_2^{-g_2}Z_1^{-g_1}$ as operators on $M_\fm$.
\end{lem}
\begin{proof}
We have
\[a_g^*a_g=(Z_1^{g_1}Z_2^{g_2})^*Z_1^{g_1}Z_2^{g_2}=Z_2^{-g_2}Z_1^{-g_1}
Z_1^{g_1}Z_2^{g_2}=Z_1^{-g_1}Z_1^{g_1}Z_2^{-g_2}Z_2^{g_2},\]
by Lemma \ref{lem:Z1kZ2lswap}. Thus by Lemma \ref{lem:ZikZiladd},
\[\varphi_\fm(a_g^*a_g)=(1-q_1)^{-|g_1|}\la^{|g_2|}q_2^{(1-2g_2+g_2)|g_2|/2}\]
which proves the formula. The last statement is immediate.
\end{proof}

We consider the three subcases corresponding to the rank of
the free abelian group $K_2$.
\subsubsection{The case $\fm=\fn_\la^{(2)}, \la\neq 0, \rank K_2=0$}
$G_\fm=K_2=\{0\}$ so $B_\fm^{(1)}=R$ which is commutative,
hence $M_\fm=\C v_0$ for some $v_0$, and $S=\Z^2$. Thus $C=\{a_gv_0\; |\; g\in\Z^2\}$
is a basis for $M$ and using Lemma \ref{lem:ZikZiladd} and Lemma
\ref{lem:Z1kZ2lswap} we obtain that the action of $X_i$ is given by
\begin{equation}\label{eq:cl_n2r0X}\begin{split}
X_1a_gv_0&=\begin{cases}
a_{g+e_1}v_0,&g_1\ge 0,\\
(1-q_1)^{-1}a_{g+e_1}v_0,&g_1<0,
\end{cases}\\
X_2a_gv_0&=\begin{cases}
(q_1\la_{12})^{-g_1}a_{g+e_2}v_0,&g_2\ge0,\\
(q_1\la_{12})^{-g_1}\la q_2^{-g_2} a_{g+e_2}v_0,&g_2<0.
\end{cases}
\end{split}\end{equation}
The action of $Y_i$ on the basis is deduced uniquely from
\begin{equation}\label{eq:cl_n2r0Y}\begin{split}
Y_1X_1a_gv_0&=(1-q_1)^{-1}a_gv_0,\\
Y_2X_2a_gv_0&=\la q_1^{-g_1}q_2^{-g_2}a_gv_0,
\end{split}\end{equation}
which hold by (\ref{eq:si1ksi2ligen}) and (\ref{eq:si1ksi2ligen2}).

\subsubsection{The case $\fm=\fn_\la^{(2)}, \la\neq 0, \rank K_2=1$}
Let $(a,b)$ be a basis element.
Since $G_\fm=K_2$ which is of rank one,
$B_\fm^{(1)}\simeq\C[T,T^{-1}]$ by Corollary \ref{cor:Bm1} so
$M_\fm$ is one-dimensional. As before
we let $M_\fm=\C v_0$. Then
$Z_1^aZ_2^bv_0=\rho v_0$ for some $\rho\in\C^*$.

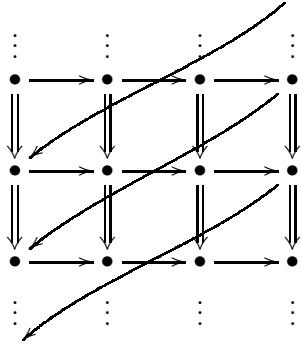
\begin{figure}
\[\xymatrix{
&&& \ar@(dl,ur)[ddlll] \\
\bullet \ar[r]  \ar@{}[u]|\vdots&
\bullet \ar[r]  \ar@{}[u]|\vdots&
\bullet \ar[r]  \ar@{}[u]|\vdots&
\bullet \ar@{}[u]|\vdots \ar@(dl,ur)[ddlll] \\
\bullet \ar[r]  \ar@{<=}[u]&
\bullet \ar[r]  \ar@{<=}[u]&
\bullet \ar[r]  \ar@{<=}[u]&
\bullet \ar@{<=}[u] \ar@(dl,ur)[ddlll]\\
\bullet \ar[r]  \ar@{<=}[u]&
\bullet \ar[r]  \ar@{<=}[u]&
\bullet \ar[r]  \ar@{<=}[u]&
\bullet \ar@{<=}[u] \\
 \ar@{}[u]|\vdots&
 \ar@{}[u]|\vdots&
 \ar@{}[u]|\vdots&
 \ar@{}[u]|\vdots
}\]
\caption{Example of a weight diagram for $M$ when $\fm=\fn^{(2)}_\la$, $\la\neq 0$ and
$\rank K_2=1$. Here $a=4$, $b=-2$.
The action of $X_1$ is indicated by $\rightarrow$ arrows, while
$\Rightarrow$ arrows are used for $X_2$.}
\label{fig:rank1}
\end{figure}

We assume $a\neq 0$. The case $b\neq 0$
can be treated similarly.
By changing basis, we can assume that $a>0$.
Choose $S=\{0,1,\ldots,a-1\}\times\Z$.
The corresponding basis for $M$ is
\[C=\{w_{ij}:=X_1^iZ_2^jv_0\; |\; 0\le i\le a-1, j\in\Z\}.\]
We now aim to apply Theorem \ref{thm:action}.
If $0\le i<a-1$ then clearly $X_1w_{ij}=w_{i+1,j}$.
And
\[X_1w_{a-1,j}=X_1^aZ_2^jv_0\in\C Z_2^{j-b}v_0=\C w_{0,j-b}.\]
We want to compute the coefficient of $w_{0,j-b}$.
Similarly to the proof of Theorem \ref{thm:action} we have,
using Lemma \ref{lem:rg}, Lemma \ref{lem:Z1kZ2lswap} and
(\ref{eq:Z2kZ2ladd}),
\begin{align*}
X_1w_{a-1,j}&=Z_1^aZ_2^jv_0=(Z_1^aZ_2^jr_{(a,b)}Z_2^{-b}Z_1^{-a})
Z_1^aZ_2^bv_0=\\
&=r_{(a,b)}(q_1\la_{12})^{ja}q_1^{a\cdot\overline{-b}}
\la_{12}^{-ab} Z_2^jZ_2^{-b} Z_1^aZ_1^{-a} \rho v_0=\\
&=(\la^{-1}q_2^{(b-1)/2})^{|b|}
q_1^{a(j+\overline{-b})}\la_{12}^{a(j-b)}\rho C_0 w_{0,j-b},
\end{align*}
where
\[C_0=
\begin{cases}
1,& b\le 0,\\
\la^{\min\{j,b\}}q_2^{(1+2b-\min\{j,b\})\min\{j,b\}/2},&b>0.
\end{cases}\]
Using Lemma \ref{lem:Z1kZ2lswap} one easily get the action of
$X_2$ on the basis. We conclude that
\begin{equation}\label{eq:cl_n2r1X}\begin{split}
X_1w_{ij}&=\begin{cases}
w_{i+1,j},&0\le i<a-1,\\
(\la^{-1}q_2^{(b-1)/2})^{|b|}q_1^{a(j+\overline{-b})}\la_{12}^{a(j-b)}\rho
C_0 w_{0,j-b},& i=a-1,
\end{cases}\\
X_2w_{ij}&=\begin{cases}
q_1^{-i}\la_{21}^i w_{i,j+1},&j\ge 0,\\
q_1^{-i}\la_{21}^i\la q_2^j w_{i,j+1},&j<0.
\end{cases}
\end{split}\end{equation}
The action of the $Y_i$ is uniquely determined by
\begin{equation}\label{eq:cl_n2r1Y}\begin{split}
Y_1X_1v_{ij}&=(1-q_1)^{-1}v_{ij},\\
Y_2X_2v_{ij}&=\la q_1^{-i}q_2^{-j} v_{ij},
\end{split}\end{equation}
which hold by (\ref{eq:si1ksi2ligen})--(\ref{eq:si1ksi2ligen2}).
See Figure \ref{fig:rank1} for a visual representation.

\subsubsection{The case $\fm=\fn_\la^{(2)}, \la\neq 0, \rank K_2=2$}
Let $s_1=\mathbf{a}=(a_1,a_2)$,
$s_2=\mathbf{b}=(b_1,b_2)$ be a basis for $G_\fm=K_2$ over $\Z$.
We can assume that $a_1,b_1\ge 0$ and that
$d:=\left|\begin{smallmatrix}a_1&b_1\\a_2&b_2\end{smallmatrix}\right| >0$.

By Corollary \ref{cor:Bm1}, $B_\fm^{(1)}\simeq T_\nu$ for some $\nu$
which we will now determine.
Using Lemma \ref{lem:Z1kZ2lswap} and Lemma \ref{lem:ZikZilswap}
we have, as operators on $M_\fm$,
\begin{align*}
Z_1^{a_1}Z_2^{a_2}Z_1^{b_1}Z_2^{b_2}&=
q_1^{-b_1\overline{a_1}}\la_{12}^{-b_1a_2}c(a_2,b_2)Z_1^{b_1}Z_1^{a_1}
Z_2^{b_2}Z_2^{a_2}=\\
&=q_1^{a_1\overline{b_2}-b_1\overline{a_2}}\la_{12}^{a_1b_2-b_1a_2}c(a_2,b_2)
Z_1^{b_1}Z_2^{b_2}
Z_1^{a_1}Z_2^{a_2}.
\end{align*}
We conclude that $B_\fm^{(1)}\simeq T_\nu$ where
\begin{equation}\label{eq:nu}
\nu = \la_{12}^d q_1^{a_1\overline{b_2}-b_1\overline{a_2}}c(a_2,b_2).
\end{equation}
The function $c$ was defined in (\ref{eq:ckl}), $d=a_1b_2-b_1a_2$ and
$\overline{k}:=\max\{0,k\}$ for $k\in\Z$.
For $M_\fm$ to be finite-dimensional it is thus necessary that this $\nu$
is a root of unity. Assume this and let $r$ denote its order. Then
$\dim M_\fm=r$. Let
\begin{equation}
\label{eq:rank2Mmbasis}
\{v_0,v_1,\ldots,v_{r-1}\}
\end{equation}
be a basis such that
\begin{align}
\label{eq:rank2localaction1}
Z_1^{a_1}Z_2^{a_2}v_j&=\nu^j\rho v_j,\\
\label{eq:rank2localaction2}
Z_1^{b_1}Z_2^{b_2}v_j&=\begin{cases}v_{j+1},&0\le j<r-1,\\
\mu v_0,&j=r-1,\end{cases}
\end{align}
where $\rho,\mu\in\C^*$.

The next step is to determine a set $S\subseteq\tilde G_\fm=\Z^2$
of representatives for the set of cosets $\tilde G_\fm / G_\fm=\Z^2/K_2$
which makes it possible to write down the action of the algebra later.
We proceed as follows.

Recall that $K_2=\Z\cdot (a_1,a_2)\oplus\Z\cdot (b_1,b_2)$.
Let $d_1$ be the smallest positive integer such that $(d_1,0)\in K_2$.
We claim that $d_1=d/\GCD(a_2,b_2)$. Indeed $d_1$ must be of the
form $ka_1+lb_1$ where $k,l\in\Z$ and $ka_2+lb_2=0$ with
$\GCD(k,l)=1$. For such $k,l$, $k|b_2$, $l|a_2$ and $b_2/k=-a_2/l=:p>0$.
Then $\GCD(a_2/p,b_2/p)=1$ which implies that $\GCD(a_2,b_2)=p$.
Thus $d_1=ka_1+lb_1=(b_2a_1-a_2b_1)/p=d/\GCD(a_2,b_2)$ as claimed.

Next, let $d_2$ denote the smallest positive integer such that
some $K_2$-translation of $(0,d_2)$ lies on the $x$-axis, i.e. such that
\[\big((0,d_2)+K_2\big)\cap \Z\times\{0\}\neq\emptyset.\]
Such an integer exists because if we write $\GCD(a_2,b_2)=ka_2+lb_2$,
then
\[(0,ka_2+lb_2)-k(a_1,a_2)-l(b_1,b_2)=(-ka_1-lb_1,0).\]
On the other hand if $(0,d_2)+k\mathbf{a}+l\mathbf{b}\in\Z\times\{0\}$,
i.e. if $d_2=ka_2+lb_2$, then $\GCD(a_2,b_2)|d_2$. Therefore
$d_2=\GCD(a_2,b_2)$.

We also see that for any point in $\Z^2$ of the form $(x,d_2)$ there
is a $g\in K_2$ such that $(x,d_2)+g\in\Z\times\{0\}$. Also,
$(d_1,0)\in K_2$ so for any point of the form $(d_1,y)$ there
is a $g\in K_2$ (namely $(-d_1,0)$) such that $(d_1,y)+g\in\{0\}\times\Z$.

Suppose now that for some $k,l\in\Z$,
\[k(a_1,a_2)+l(b_1,b_2)\in K_2\cap\{0,1,\ldots,d_1-1\}\times
\{0,1,\ldots,d_2-1\}.\]
Then we would have
$(0,ka_2+lb_2)-(k\mathbf{a}+l\mathbf{b})\in\Z\times\{0\}$
and $ka_2+lb_2\in\{0,1,\ldots,d_2-1\}$ which contradicts the minimality
of $d_2$ unless $ka_2+lb_2=0$. But in this case $(ka_1+lb_1,0)\in K_2$
which contradicts the minimality of $d_1$ unless $ka_1+lb_1=0$. Hence
$K_2\cap\{0,1,\ldots,d_1-1\}\times\{0,1,\ldots,d_2-1\}=\{(0,0)\}$.
We have shown that
\[S:=\{0,1,\ldots,d_1-1\}\times\{0,1,\ldots,d_2-1\}\]
is a set of representatives for $\Z^2 / K_2$.
In particular we get from Corollary \ref{cor:dimM} that $\dim M$ is
finite and
\[\dim M/\dim M_\fm=|S|=d_1d_2=a_1b_2-b_1a_2.\]

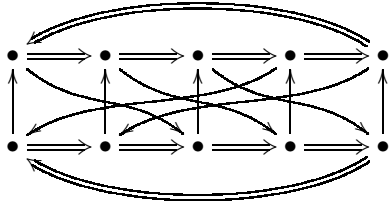
\begin{figure}
\[\xymatrix{
\bullet \ar@{=>}[r] \ar@(dr,ul)[drr]&
\bullet \ar@{=>}[r] \ar@(dr,ul)[drr]&
\bullet \ar@{=>}[r] \ar@(dr,ul)[drr]&
\bullet \ar@{=>}[r] \ar@(dl,ur)[dlll]&
\bullet \ar@(ul,ur)@{=>}[llll] \ar@(dl,ur)[dlll] \\
\bullet \ar@{=>}[r] \ar[u] &
\bullet \ar@{=>}[r] \ar[u] &
\bullet \ar@{=>}[r] \ar[u] &
\bullet \ar@{=>}[r] \ar[u] &
\bullet \ar@(dl,dr)@{=>}[llll] \ar[u]
}\]
\caption{An example of the action on $\supp(M)$
when $\fm=\fn^{(2)}_\la$, $\la\neq 0$ and $\rank K_2=2$. Here
$\mathbf{a}=(2,-2)$, $\mathbf{b}=(3,2)$, $d_1=5, d_2=2$ and $s=2$.
The $\Rightarrow$ arrows indicate the action of $X_1$ and
the $\rightarrow$ arrows show the action of $X_2$.}
\label{fig:rank2}
\end{figure}

We fix now integers $a_2',b_2'$ such that
\begin{equation}\label{eq:d2}
d_2=\GCD(a_2,b_2)=a_2'a_2+b_2'b_2
\end{equation}
and such that
$-a_2'a_1-b_2'b_1\in\{0,1,\ldots,d_1-1\}$. This can be done because
for any $p\in\Z$, $(a_2'',b_2''):=
(a_2'+pb_2/d_2,b_2'-pa_2/d_2)$ also satisfies $a_2''a_2+b_2''b_2=d_2$
but now
\[-a_2''a_1-b_2''b_1=-(a_2'+pb_2/d_2)a_1-(b_2'-pa_2/d_2)b_1=
-a_2'a_1-b_2'b_1-pd_1.\]
We set
\begin{equation}\label{eq:s}
s=-a_2'a_1-b_2'b_1.
\end{equation}
Let $(i,j)\in S$.
We have the following reductions in $\Z^2$ modulo $K_2$.
\begin{align*}
(1,0)+(i,j)&=\begin{cases}
(i+1,j),& 0\le i<d_1-1,\\
(0,j),&i=d_1-1,\end{cases}\\
(0,1)+(i,j)&=\begin{cases}
(i,j+1),&0\le j<d_2-1,\\
(i+s,0),&j=d_2-1, i+s\le d_1-1,\\
(i+s-d_1,0),&j=d_2-1, j+s>d_1-1.
\end{cases}
\end{align*}
From this we can understand how the $X_i$ act on the support of $M$,
see Figure \ref{fig:rank2} for an example.
By Theorem \ref{thm:basis} the set
\[C=\{w_{ijk}:=X_1^iX_2^jv_k\; |\; 0\le i<d_1, 0\le j<d_2, 0\le k<r\}\]
is a basis for $M$ where $v_k$ is the basis (\ref{eq:rank2Mmbasis})
for $M_\fm$.

If $0\le i<d_1-1$ we clearly have $X_1w_{ijk}=w_{i+1,j,k}$.
Suppose $i=d_1-1$. Then by Lemma \ref{lem:Z1kZ2lswap},
\[X_1w_{ijk}=X_1^{d_1}X_2^jv_k=q_1^{d_1j}\la_{12}^{d_1j}X_2^jX_1^{d_1}v_k.\]
Thus we must express $X_1^{d_1}$ in terms of $Z_1^{a_1}Z_2^{a_2}$
and $Z_1^{b_1}Z_2^{b_2}$. Since $(d_1,0)=b_2/d_2\mathbf{a}-a_2/d_2
\mathbf{b}$ we have
\begin{equation}\label{eq:C1inverse}
(Z_1^{a_1}Z_2^{a_2})^{b_2/d_2}(Z_1^{b_1}Z_2^{b_2})^{-a_2/d_2}=
C_1^{-1}X_1^{d_1}
\end{equation}
as operators on $M_\fm$ for some constant $C_1^{-1}$ which we must calculate.

\begin{lem} \label{lem:C1}
The constant $C_1$ defined in (\ref{eq:C1inverse})
is given by
\begin{multline}\label{eq:C1inverse2}
C_1^{-1}=r_\mathbf{a}^{\overline{-b_2/d_2}}(q_1^{-a_1\overline{a_2}}
\la_{12}^{-a_1a_2})^{\frac{b_2}{d_2}(\frac{b_2}{d_2}-1)/2}\cdot
r_\mathbf{b}^{\overline{a_2/d_2}}(q_1^{-b_1\overline{b_2}}
\la_{12}^{-b_1b_2})^{\frac{a_2}{d_2}(\frac{a_2}{d_2}+1)/2}\cdot\\ \cdot
q_1^{b_1a_2\overline{a_2b_2}/d_2^2}\la_{12}^{b_1a_2^2b_2/d_2^2}
r_{(0,-b_2a_2/d_2)}^{-1}C_1',
\end{multline}
where the $r_g$, $g\in \Z^2$ are given by (\ref{eq:rg}),
\[C_1'=\begin{cases}
(1-q_1)^{-\min\{|a_1b_2/d_2|,|b_1a_2/d_2|\}},&a_2b_2>0,\\
1,&a_2b_2\le 0,
\end{cases}\]
$\overline{k}=\max\{0,k\}$ for $k\in\Z$ and $d_2=\GCD(a_2,b_2)$.
\end{lem}
\begin{proof}
If $b_2\ge 0$ for example, we have by Lemma \ref{lem:Z1kZ2lswap}
\begin{multline*}
(Z_1^{a_1}Z_2^{a_2})^{b_2/d_2}=q_1^{-a_1\overline{a_2}}\la_{12}^{-a_1a_2}
\cdot(q_1^{-a_1\overline{a_2}}\la_{12}^{-a_1a_2})^2\cdot\ldots\\
\ldots\cdot
(q_1^{-a_1\overline{a_2}}\la_{12}^{-a_1a_2})^{b_2/d_2-1}Z_1^{a_1b_2/d_2}
Z_2^{a_2b_2/d_2}=\\
=(q_1^{-a_1\overline{a_2}}\la_{12}^{-a_1a_2})^{\frac{b_2}{d_2}(\frac{b_2}{d_2}-1)/2}
Z_1^{a_1b_2/d_2}Z_2^{a_2b_2/d_2}
\end{multline*}
When $b_2<0$ we get a similar calculation where $r_\mathbf{a}^{-b_2/d_2}$
appears by Lemma \ref{lem:rg}. $(Z_1^{b_1}Z_2^{b_2})^{-a_2/d_2}$ can
analogously be expressed as a multiple of $Z_1^{-b_1a_2/d_2}Z_2^{-b_2a_2/d_2}$.
We then commute $Z_2^{a_2b_2/d_2}$ and $Z_1^{-b_1a_2/d_2}$ using
Lemma \ref{lem:Z1kZ2lswap}. As a last step we use Lemma \ref{lem:ZikZiladd}
and obtain two more factors.
\end{proof}

We conclude that
\begin{align*}
X_1w_{ijk}&=\begin{cases}
w_{i+1,j,k},&i<d_1-1,\\
q_1^{jd_1}\la_{12}^{jd_2}C_1\nu^{b_2/d_2 k_1''}\rho^{b_2/d_2}
\mu^{k_1'} w_{0,j,k_1''},& i=d_1-1.
\end{cases}
\end{align*}
Here
\begin{equation}\label{eq:k1}
k-a_2/d_2=rk_1'+k_1''\quad \text{with } 0\le k_1''<r.
\end{equation}

Next we turn to the description of how $X_2$ acts on the basis $C$.
If $0\le j<d_2-1$ we have $X_2w_{ijk}=q_1^{-i}\la_{12}^{-i}w_{i,j+1,k}$
by Lemma \ref{lem:Z1kZ2lswap}.
Suppose $j=d_2-1$. Then, as in the first step of the proof
of Theorem \ref{thm:action},
\begin{equation}\label{eq:X2d2}
X_2w_{ijk}=q_1^{-i}\la_{12}^{-i}X_1^{i}X_2^{d_2}v_k=
q_1^{-i}\la_{12}^{-i}X_1^{i}
(X_2^{d_2}r_{(-s,d_2)}Z_2^{-d_2}Z_1^s)(Z_1^{-s}Z_2^{d_2})v_k.
\end{equation}
By (\ref{eq:Z2kZ2ladd}) and (\ref{eq:rg}),
\begin{align}\label{eq:rank2gen0}
X_2^{d_2}r_{(-s,d_2)}Z_2^{-d_2}Z_1^s&=
r_{(-s,d_2)}r_{(0,-d_2)}^{-1}Z_1^s=\\ \nonumber
&=(1-q_1)^s(\la^{-1}q_2^{(d_2-1)/2})^{d_2}(\la^{-1}q_2^{(-d_2-1)/2})^{d_2}Z_1^s=\\
&=(1-q_1)^s(\la^2q_2)^{-d_2}Z_1^s. \nonumber
\end{align}
We must express $Z_1^{-s}Z_2^{d_2}$ in the generators of the algebra
$B_\fm^{(1)}$ in order to calculate its action on $v_k$.
\begin{equation}\label{eq:rank2gen1}
(Z_1^{a_1}Z_2^{a_2})^{a_2'}(Z_1^{b_1}Z_2^{b_2})^{b_2'}
= C_2^{-1} Z_1^{-s}Z_2^{d_2},
\end{equation}
for some $C_2\in\C^*$ since the degree on both sides are equal by
(\ref{eq:d2}) and (\ref{eq:s}).
Similarly to the proof of Lemma \ref{lem:C1},
\begin{multline}\label{eq:C2inverse}
C_2^{-1}=r_\mathbf{a}^{\overline{-a_2'}}(q_1^{-a_1\overline{a_2}}
\la_{12}^{-a_1a_2})^{a_2'(a_2'-1)/2}\cdot
r_\mathbf{b}^{\overline{-b_2'}}(q_1^{-b_1\overline{b_2}}
\la_{12}^{-b_1b_2})^{b_2'(b_2'-1)/2}\cdot\\
\cdot q_1^{-b_1b_2'\overline{a_2a_2'}}\la^{-b_1b_2'a_2a_2'}C_2'C_2'',
\end{multline}
and
\begin{align*}
C_2'&=\begin{cases}
1,&a_2'b_2'\ge 0,\\
(1-q_1)^{-\min\{|a_1a_2'|,|b_1b_2'|\}},&a_2'b_2'<0,
\end{cases}\\
C_2''&=\begin{cases}
1,&a_2a_2'b_2b_2'\ge 0,\\
\la^{m'}q_2^{(1-2b_2b_2'+(\sgn b_2b_2')m')m'/2},& a_2a_2'b_2b_2'<0,
\end{cases}
\end{align*}
where $m'=\min\{|a_2a_2'|,|b_2b_2'|\}$.
Furthermore, letting
\begin{equation}\label{eq:k2}
b_2'+k=rk_2'+k_2'',\quad\text{where }0\le k_2''<r
\end{equation}
we have by (\ref{eq:rank2localaction1})--(\ref{eq:rank2localaction2}),
\begin{equation}\label{eq:rank2gen2}
(Z_1^{a_1}Z_2^{a_2})^{a_2'}(Z_1^{b_1}Z_2^{b_2})^{b_2'}v_k=
\nu^{a_2'k_2''}\rho^{a_2'}\mu^{k_2'}v_{k_2''}.
\end{equation}

If $i+s\le d_1-1$ we can now write down the action of $X_2$ on $w_{ijk}$
by combining (\ref{eq:X2d2})--(\ref{eq:C2inverse}), (\ref{eq:rank2gen2})
to get a multiple of $w_{i+s,0,k_2''}$. However if $i+s>d_1-1$,
we must reduce further because then $(i+s,0)\notin S$. Let
\begin{equation}\label{eq:k3}
k_2''-a_2/d_2=rk_3'+k_3'',\quad\text{where }0\le k_3''<r.
\end{equation}
Then by the calculations for the action of $X_1^{d_1}$ on $M_\fm$,
\[X_1^{d_1}v_{k_2''}=X_1^{i+s-d_1}X_1^{d_1}v_{k_2''}=
C_1\mu^{k_3'}\nu^{k_3''b_2/d_2}\rho^{b_2/d_2}w_{i+s-d_1,0,k_3''}.\]

Summing up, $M$ has a basis
\[\{w_{ijk}\; |\; 0\le i<d_1, 0\le j<d_2, 0\le k<r\}\]
and $X_1, X_2$ act on this basis as follows.
\begin{equation}\label{eq:cl_n2r2X}\begin{split}
X_1w_{ijk}&=\begin{cases}
w_{i+1,j,k},&i<d_1-1,\\
q_1^{jd_1}\la_{12}^{jd_2}C_1\nu^{b_2/d_2 k_1''}\rho^{b_2/d_2}
\mu^{k_1'} w_{0,j,k_1''},& i=d_1-1.
\end{cases}\\
X_2w_{ijk}&=(q_1\la_{12})^{-i}\cdot\\&\cdot
\begin{cases}
w_{i,j+1,k},&\\ \qquad \text{if }0\le j<d_2-1, &\\
(1-q_1)^s(\la^2q_2)^{-d_2}C_2
\nu^{a_2'k_2''}\rho^{a_2'}\mu^{k_2'}
w_{i+s,0,k_2''},&\\ \qquad \text{if } j=d_2-1 \text{ and } i+s\le d_1-1,&\\
(1-q_1)^s(\la^2q_2)^{-d_2}C_2
\nu^{a_2'k_2''+k_3''b_2/d_2}\rho^{a_2'+b_2/d_2}\mu^{k_2'+k_3'}
C_1w_{i+s-d_1,0,k_3''}
,&\\ \qquad \text{if }j=d_2-1 \text{ and } i+s>d_1-1, &
\end{cases}
\end{split}\end{equation}
where $C_1$ is given by (\ref{eq:C1inverse2}), $C_2$ by
(\ref{eq:C2inverse}) and $\nu$ by (\ref{eq:nu}). The parameters
$\rho$ and $\mu$ comes from the action (\ref{eq:rank2localaction1}),
(\ref{eq:rank2localaction2}) of $B_\fm^{(1)}$ on $M_\fm$ and
$k_i', k_i''$ are defined in (\ref{eq:k1}), (\ref{eq:k2}) and (\ref{eq:k3}).

The action of the $Y_i$ is uniquely determined by
\begin{equation}\label{eq:cl_n2r2Y}\begin{split}
Y_1X_1w_{ijk}&=(1-q_1)^{-1}w_{ijk},\\
Y_2X_2w_{ijk}&=\la q_1^{-i}q_2^{-j} w_{ijk}.
\end{split}\end{equation}

We remark that the case $q_1=q_2$ corresponds to $\mathbf{a}=(a_1,a_2)=(1,-1)$.
Then $d_2=1$, $d_1=d=|b_1+b_2|$ and $s=1$. $X_1$ and $X_2$ will act on the
support in the same direction, cyclically as in Figure \ref{fig:rank2cyc}.
\begin{figure}
\[\xymatrix{
\bullet \ar[r]&
\bullet \ar[r]&
\bullet \ar@{}[r]|{\ldots\ldots}&
\bullet \ar[r]&
\bullet \ar@(ul,ur)[llll]
}\]
\caption{Weight diagram when $\fm=\fn^{(2)}_\la$, $\la\neq 0$, $\rank K_2=2$
and $q_1=q_2$.}
\label{fig:rank2cyc}
\end{figure}
The explicit action can be deduced from the above more general case
noting that here $k_2''=k$, $k_2'=0$ and
\begin{align*}
k_1'=k_3'&=\begin{cases}
0,&k<r-1,\\ 1,&k=r-1,
\end{cases}&
k_1''=k_3''&=\begin{cases}
k,&k<r-1,\\ 0,&k=r-1.
\end{cases}
\end{align*}

%
%
%   F o u r t h   c a s e
%
%
\subsection{The case
 $\fm \not\in \{\fn_{\mu}^{(i)}\;|\; \mu\in\C, i=1,2\}$}
This is the generic case. We have $\Z^2_\fm=Q$
by Corollary \ref{cor:Z2fm}.
Our statements here generalize without any problem
to the case of arbitrary rank.

Assume first that the $q_i$ are roots of unity of
orders $o_i$ ($i=1,2$) and that
$\omega$ does not contain any $1$-breaks
or $2$-breaks. Then by Corollary \ref{cor:breaks}
and Proposition \ref{prop:tildeGm} we have
$\tilde G_\fm=\Z^2$. Thus $G_\fm=(o_1\Z)\times(o_2\Z)$.
Moreover,
\[X_1^{o_1}X_2^{o_2}=\la_{12}^{o_1o_2}X_2^{o_2}X_1^{o_1}\]
so $B_\fm^{(1)}\simeq T_{\la_{12}^{o_1o_2}}$
by Corollary \ref{cor:Bm1}. This algebra
has only finite dimensional representations if $\la_{12}^{o_1o_2}$ is
a root of unity. Assuming this, let $r$ be the order
of $\la_{12}^{o_1o_2}$. Then there are $\rho,\mu\in\C^*$ and
$M_\fm$ has a basis $v_0,v_1,\ldots,v_{r-1}$ such that
\begin{align*}
X_1^{o_1}v_i&=\la_{12}^{io_1o_2}\rho v_i\\
X_2^{o_1}v_i&=\begin{cases} v_{i+1}&0\le i<p-1\\ \mu v_0&i=p-1
\end{cases}
\end{align*}
Choose $S=\{0,1,\ldots,o_1-1\}\times\{0,1,\ldots,o_2-1\}$.
The corresponding basis for $M$ is
$C=\{w_{ijk}:=X_1^iX_2^jv_k\; |\; 0\le i<o_1, 0\le j<o_2, 0\le k<r\}$.
The following formulas are easily deduced using
(\ref{eq:qweylrel1})--(\ref{eq:qweylrel3}).
\begin{equation}\label{eq:cl_giX}\begin{split}
X_1 w_{ijk}&=
\begin{cases}
w_{i+1,j,k},&k<o_1-1,\\
\la_{12}^{o_1(o_2k+j)}\rho w_{0jk},& k=o_1-1,
\end{cases}\\
X_2 w_{ijk}&= (q_1\la_{12})^{-i}\cdot
\begin{cases}
 w_{i,j+1,l},& l<o_2-1,\\
 w_{i,0,l+1},& l=o_2-1, i<r-1,\\
 \mu w_{i00},& l=o_2-1, i=r-1.
\end{cases}
\end{split}\end{equation}
The action of $Y_1,Y_2$ is determined by
\begin{equation}\label{eq:cl_giY}\begin{split}
Y_1X_1w_{ijk}&=q_1^{-i}(\al_1-[i]_{q_1})w_{ijk},\\
Y_2X_2w_{ijk}&=q_1^{-i}q_2^{-j}(\al_2-[j]_{q_2}(1+(q_1-1)\al_1))w_{ijk}.
\end{split}\end{equation}

In all other cases one can show using the same argument that
$\dim M_\fn=1$ for all $\fn\in\supp(M)$ and that $M$ can be realized in a
vector space with basis $\{w_{ij}\}_{(i,j)\in I}$, where $I=I_1\times I_2$
is one of the following sets
\begin{align*}
&\N_{d_1}\times \N_{d_2},\quad \N_{d_1}\times \Z^\pm,\quad \Z^\pm\times\N_{d_2},\quad \Z\times\Z,\\
&\Z^\pm\times\Z,\quad \Z\times\Z^\pm,\quad \Z^\pm\times\Z^\pm, \quad\Z^\pm\times\Z^\mp,
\end{align*}
where $\N_d=\{0,1,\ldots,d-1\}$, $\Z^\pm=\{k\in\Z\; |\; \pm k\ge 0\}$
and $d_i$ is the order of $q_i$ if finite.
The action of the generators is given by the following formulas.
\begin{equation}\label{eq:cl_giiX}\begin{split}
X_1w_{ij}&=\begin{cases}
w_{i+1,j},&(i+1,j)\in I,\\
\rho\la_{12}^{d_1j}w_{0,j},& (i+1,j)\notin I, I_1=\N_{d_1} \text{ and }
\al_1\neq [i]_{q_1},\\
0,&\text{otherwise},
\end{cases}\\
X_2w_{ij}&=(q_1\la_{12})^{-i}\cdot\begin{cases}
w_{i,j+1},&(i,j+1)\in I,\\
\mu w_{i,0},& (i,j+1)\notin I, I_2=\N_{d_2}\\&\quad \text{and }
\al_2\neq [j]_{q_2}(1+(q_1-1)\al_1),\\
0,&\text{otherwise},
\end{cases}
\end{split}\end{equation}
\begin{equation}\label{eq:cl_giiY}\begin{split}
Y_1w_{ij}&=q_1^{-i+1}(\al_1-[i-1]_{q_1})\cdot\\ &\cdot \begin{cases}
w_{i-1,j},&(i-1,j)\in I,\\
(\rho \la_{12}^{d_1j})^{-1}w_{d_1-1,j},&(i-1,j)\notin I, I_1=\N_{d_1}
\text{ and } \al_1\neq [i-1]_{q_1},\\
0,&\text{otherwise},
\end{cases}\\
Y_2w_{ij}&=\la_{12}^{-i}q_2^{-j+1}(\al_2-[j-1]_{q_2}(1+(q_1-1)\al_1))\cdot\\
&\cdot
\begin{cases}
w_{i,j+1},&(i,j+1)\in I,\\
\mu^{-1}w_{i,d_2-1},&(i,j+1)\notin I, I_1=\N_{d_2} \\&\quad\text{ and }
\al_2\neq [j-1]_{q_2}(1+(q_1-1)\al_1),\\
0,&\text{otherwise}.
\end{cases}
\end{split}\end{equation}

Thus we have proved the following result.
\begin{thm}\label{thm:QWAcl}
Let $A$ be a quantized Weyl algebra
of rank two with arbitrary parameters $q_1, q_2\in\C\backslash\{0,1\}$.
Then any simple weight $A$-module with no proper inner breaks
is isomorphic to one of the modules
defined by formulas (\ref{eq:cl_n0}), (\ref{eq:cl_n1i}),
(\ref{eq:cl_n1iihi}), (\ref{eq:cl_n1iilo}), (\ref{eq:cl_n1iii}),
(\ref{eq:cl_n1iv}), (\ref{eq:cl_n2r0X}-\ref{eq:cl_n2r0Y}),
(\ref{eq:cl_n2r1X}-\ref{eq:cl_n2r1Y}),
(\ref{eq:cl_n2r2X}-\ref{eq:cl_n2r2Y}),
(\ref{eq:cl_giX}-\ref{eq:cl_giY}) or
(\ref{eq:cl_giiX}-\ref{eq:cl_giiY}).
\end{thm}

\section*{Acknowledgments}
The author is very grateful to his supervisor L. Turowska
for introducing him to the subject and for many stimulating conversations.
Thanks are also
due to D. Proskurin for interesting discussions and to V. Mazorchuk for helpful
comments.


\begin{thebibliography}{99}
\bibitem{Ba} Bavula, V., \emph{Finite-dimensionality of ${\rm Ext}\sp n$ and
${\rm Tor}\sb n$ of simple modules over a class of algebras.} (Russian)
Funktsional. Anal. i Prilozhen.  25  (1991), no. 3, 80--82;  translation in
Funct. Anal. Appl.  25  (1991), no. 3, 229--230 (1992)
\bibitem{B} Bavula, V., \emph{Generalized Weyl algebras and their representations},
St. Petersb. Math. J., 4, (1993), 71--93.
\bibitem{B2} Bavula, V., Bekkert, V., \emph{Indecomposable representations of generalized
Weyl algebras}, Comm. Algebra 28 (2000), 5067-5100.
\bibitem{BB} Brenken, B., \emph{A classification of some noncommutative tori},
The Rocky Mountain Journal of Mathematics, Vol. 20, No. 2 (1990), 389--397.
\bibitem{DFO} Drozd, Yu., Futorny, V., Ovsienko, S.,
\emph{Harish-Chandra subalgebras and Gelfand-Zetlin modules},
in: \emph{Finite-dimensional algebras and related topics},
(Ottawa, ON, 1992),
%Kluwer Academic, Dordrecht, 1994.
NATO Adv. Sci. Inst. Ser. C. Math. Phys. Sci. 424, (1994), 79--93.
\bibitem{DGO} Drozd, Yu., Guzner, B., Ovsienko, S.,
\emph{Weight modules over generalized Weyl algebras},
J. Algebra 184, (1996), 491--504.
\bibitem{MPT} Mazorchuk, V., Ponomarenko, M., Turowska, L., \emph{Some associative
algebras related to $U(\mathfrak{g})$ and twisted generalized Weyl algebras},
Math. Scand. 92 (2003), 5--30.
\bibitem{MT} Mazorchuk, V., Turowska, L., \emph{Simple weight modules
over twisted generalized Weyl algebras},
Comm. Alg. 27(6) (1999), 2613--2625.
\bibitem{newman} Newman, M., \emph{Integral Matrices}, Academic Press, New York, 1972.
%\bibitem{S} Smith, S.,
%\emph{A class of algebras similar to the enveloping algebra of $sl(2)$},
%Trans. AMS, 322, (1990), 285--314.
\end{thebibliography}
\end{document}